\documentclass[12pt,a4paper,final,english]{article}

\usepackage[textwidth=16cm,top=2.2cm,textheight=24.24cm]{geometry}
\usepackage[final]{pdfpages}

\usepackage[utf8]{inputenc}  
\usepackage{epic,eepic}
\usepackage{mathrsfs}   
\usepackage{eucal} 
\usepackage{xcolor,scrtime,graphicx,hyperref}
\usepackage{amsmath,amsfonts,amsmath,amssymb}
\usepackage{hyperref}
\usepackage{tikz} 
\numberwithin{equation}{section}
\numberwithin{figure}{section}
\usepackage{bef_alex}

\def\Bigset#1#2{\Big\{\; #1 \; \Big| \; #2 \; \Big\} }
\def\bigset#1#2{\big\{\; #1 \; \big| \; #2 \; \big\} }
\newcommand{\mafo}{\mathrm}

\newcommand{\Div}{\mathop{\mafo{div}}}
\newcommand{\PROB}{\mathop{\mafo{Prob}}}
\newcommand{\weak}{\rightharpoonup}
\newcommand{\weaks}{\overset{*}{\weak}}
\newcommand{\EDPweak}{\overset{\mafo{EDP}}{\weak}}
\newcommand{\Gto}{\overset{\Gamma}{\rightarrow}}
\newcommand{\Gweak}{\overset{\Gamma}{\rightharpoonup}}
\newcommand{\Gweaks}{\overset{\Gamma}{\rightharpoonup}\!\!{}^*\,}
\newcommand{\diag}{\mathop{\mafo{diag}}}
\newcommand{\CCC}{{\mathscr C}}
\newcommand{\pEweak}{\overset{\mathrm{pE}}{\weak}}

\usepackage{psfrag,bm} 

\usepackage[notref,notcite]{showkeys}

\begin{document}
\title{On microscopic origins of \\ generalized gradient
  structures\thanks{Partially supported by Einstein Stiftung Berlin, 
   ERC AdG\,267802, and DFG via SFB\,1114}}
\author{Matthias Liero\thanks{Weierstra\ss-Institut f\"ur
    Angewandte Analysis und Stochastik, Berlin},\quad Alexander
    Mielke$,{\!\!}^\dagger$ \\  Mark A. Peletier\thanks{Centre for Analysis,
    Scientific Computing and Applications, Technische Universiteit
    Eindhoven},\quad and\quad D. R. Michiel Renger$^\dagger$}
\date{\today} 

\maketitle

\begin{abstract}\noindent
%
  Classical gradient systems have a linear relation between rates and
  driving forces.  In generalized gradient systems we allow for
  arbitrary relations derived from general non-quadratic dissipation
  potentials.  This paper describes two natural origins for these
  structures.

A first microscopic origin of generalized gradient structures is given by
the theory of large-deviation principles. While Markovian diffusion
processes lead to classical gradient structures, Poissonian jump processes
give rise to cosh-type dissipation potentials. 
  
A second origin arises via a new form of con\-ver\-gence, that we call
EDP-con\-ver\-gence. Even when starting with classical gradient
systems, where the dissipation potential is a quadratic functional of
the rate, we may obtain a generalized gradient system in the
evolutionary $\Gamma$-limit.  As examples we treat (i) the limit of a
diffusion equation having a thin layer of low diffusivity, which leads
to a membrane model, and (ii) the limit of diffusion over a high
barrier, which gives a reaction-diffusion system.
\end{abstract}


\section{Introduction}
\label{s:Intro}

We consider evolution equations $\dot u = V(t,u)$ 
that are generated by gradient systems (GS). By a GS we understand a triple
$(\bfX,\calE,\calR)$, where the state space $\bfX$ is a weakly closed 
convex subset of a Banach space  containing the
states $u(t)$. The functional $\calE:\bfX\to \R\cup\{\infty\}$ is called
energy, but in applications it may be a free energy, a relative entropy, or
the negative of the entropy. Finally $\calR$ is the dissipation
potential depending on the state $u$ and the rate $\dot u$ such that
$\rmD_{\dot u} \calR(u,\dot u)\in \bfX^*$ denotes the dissipation
force. The induced evolution equation is the force balance
 \begin{equation}
  \label{eq:I.FB}
0= \rmD_{\dot u}\calR(u(t),\dot u(t)) + \rmD_u\calE(t,u(t)),
\end{equation}
where the symbol $\rmD$ denotes the (partial) Gateaux derivative or the
convex subdifferential. Quite often, we will use the dual dissipation
potential $\calR^*$ that is defined by the Legendre-Fenchel transform of
$\calR(u,\cdot)$. Then, the evolution equation can be rewritten as 
\begin{equation}
  \label{eq:I.GF}
  \dot u(t) = \rmD_\xi \calR^*\big(u(t),{-}\rmD\calE(u(t))\big),
\end{equation}
see Section \ref{ss:VarPrinc} for the details. Since $\calR$ and
$\calR^*$ are in one-to-one correspondence, we will sometimes denote 
$(\bfX,\calE,\calR)$ also by $(\bfX,\calE,\calR^*)$, in
particular if $\calR^*$ is given explicitly.

A third equivalent formulation of the gradient flow 
is given via the \emph{energy-dissipation principle (EDP)},
also called De Giorgi's $(\calR,\calR^*)$ principle, cf.\
\cite{DeMaTo80PEMS,AmGiSa05GFMS}. This states that, under suitable technical
assumptions, a curve $u:[0,T]\to \bfX$ is a solution of \eqref{eq:I.FB} or
\eqref{eq:I.GF} if and only if it satisfies the energy-dissipation
estimate
\begin{equation}
  \label{eq:I.EDE}
  \calE(u(T))+{\mathscr D}(u)\leq
\calE(u(0)) \ \text{with } {\mathscr D}(u):=\int_0^T 
\!\!\calR(u(t),\dot u(t)) + \calR^*\big(u(t), {-}\rmD
\calE(u(t)) \big) \dd t. 
\end{equation}
We call $\mathscr D$ the \emph{De Giorgi dissipation functional}. 

A GS is called \emph{classical} if the dissipation potential
$\calR(u,\cdot)$ is quadratic, i.e.\ $\calR(u,\dot u)=\frac12\langle
\bbG(u)\dot u,\dot u\rangle$ for a linear, symmetric, and positive definite
 operator $\bbG(u):\bfX\to\bfX^*$. 
If we want to emphasize that a GS is not
classical, we call it a \emph{generalized GS}. The aim of this work is
to show that generalized GS arise in two natural ways. First, it is
shown in \cite{MiPeRe14RGFL,MPPR15?LDPC} that they appear via
large-deviation principles from a microscopic $N$-particle system for
$N\to \infty$, see Section \ref{ss:Markov} for a brief summary of the
main result.  Second, generalized GS occur as suitable multiscale
limits of classical GS. \medskip

Obviously, every GS generates exactly one gradient-flow evolution equation by
\eqref{eq:I.FB} or \eqref{eq:I.GF}, but a given evolution equation
$\dot u=V(t,u)$ may be generated by many GS. If there exists at least
one such GS, we say that the evolution equation has a gradient
structure, if we do not want to specify the particular GS.
As an elementary example we treat the scalar ODE 
\[
\dot p = 1 - 2 p \ \text{ with }p(t)\in [0,1]=\PROB(\{1,2\}),
\]
which we interpret as the Kolmogorov forward equation for a Markov
process $(X_t)_{t\geq0}$ with $X_t\in \{1,2\}$. Obviously, this ODE is
generated by the GS $([0,1],\calE_2,\calR_2)$ with 
\[
\calE_2(p)= a\,\big(p-1/2\big)^2 \quad \text{and} \quad
\calR_2(p,\dot p)= \frac a2\,\dot p{}^2
\]
for any $a>0$. Of course, GS that simply differ by a scaling constant
such as $a>0$ are not considered as different. Motivated by a Markovian 
large-deviation principle, a truly different GS is
obtained for $a>0$ by
\[
\calE_\text{Mv}(p)= a\Big(p \log p +(1{-}p)\log(1{-}p)\Big)
\text{ and } 
\calR^*_\text{Mv}(p,\xi) = a\sqrt{p(1{-}p)} \,\CCC^*(\xi/a),
\]
where the function $\CCC$ and its Legendre dual $\CCC^*$ are given by 
\begin{equation}
  \label{eq:SH}
  \CCC(v) = 2 v \mathop{\mathrm{arsinh}}(v/2) - 2\sqrt{4{+}v^2} +4 \quad 
\text{and} \quad \CCC^*(\xi)=4\big(\cosh(\xi/2)-1\big). 
\end{equation}
The functions $\CCC$ and $\CCC^*$ will play a fundamental role, so we give
some elementary relations:
\begin{align*}
  &\CCC(v)=\tfrac12 v^2 + O(v^4), \quad \CCC^*(\xi)=\tfrac12 \xi^2 + 
    O(\xi^4), \quad \CCC'(v)=2\mathop{\mathrm{arsinh}}(v/2), \\[0.5em]
  & \sqrt{pq}\,\CCC^*(\log p{-}\log q)= 2\big( \sqrt p - \sqrt
  q\big)^2, \qquad \sqrt{pq}\:(\CCC^*)'(\log p{-}\log q \big) = p - q.
\end{align*} 
Indeed, using the last relation and
$\rmD\calE_\text{Mv}(p)=a\big(\log p- \log(1{-}p)\big)$ we easily find 
\[
\dot p = \rmD_\xi \calR_\text{Mv}\big(p,{-}\rmD\calE_\text{Mv}(p)\big) =
1-2p. 
\] 
Moreover, using $(\CCC^*)'(\xi)=2\sinh(\xi/2)$ we see that the evolution
takes the exponential form 
\[
\dot p= -2\sqrt{p(1{-}p)}\,\sinh\Big( \frac12 \rmD\calE_\text{Mv}(p)\Big). 
\]
This form is derived and extensively studied in \cite{BonPel14?QRIL},
it occurs in mechanics \cite[Eqn.\,(5)]{RoRaGo00WHHA} and in
chemistry, see the discussion at the end of Section \ref{sss:nonlRDS}.

The usage of generalized GS is common in the modeling of materials,
e.g.\ for plasticity, ferromagnetism, etc., where the nonsmoothness and
nonlinearity of the constitutive law $\dot u \mapsto \rmD_{\dot u}
\calR(u,\dot u)$ for the dissipative forces is essential, see Section
\ref{sss:DissMater} and the survey \cite{Miel15VAMD}. The mathematical
usage of generalized GS in smooth models such as reaction-diffusion
equations and systems is rather new. One of the remarkable origins of
gradient structures arises from the interpretation of a macroscopic
system as a Kolmogorov forward equation 
\begin{equation}
  \label{eq:I.KFe}
  \dot \rho = \bbQ^* \rho, \quad \text{where }\rho(t) \in \PROB(S),
\end{equation}
for a Markov process $(X(t))_{t\geq 0}$ on the set $S$  with
generator $\bbQ$. Considering $N$ independent particles $X_j(t)$,
$j=1,\ldots,N$ one can define the empirical process $\rho^N(t)=\frac1N
\sum_{j=1}^N \delta_{X_j(t)}\in \PROB(S)$. For $N\to \infty$ the
process $\rho^N$ converges to a solution $\rho$ of \eqref{eq:I.KFe}.
Moreover, according to the program in \cite{ADPZ11LDPW,ADPZ13LDGF} 
$\rho^N$ satisfies a large-deviation principle that gives rise to a 
rate functional 
$  \mathscr I(\rho(\cdot)) = \int_0^T \calL(\rho(t),\dot\rho(t)) \dd t$,
where $\calL$ can be characterized explicitly by $\bbQ$. The main
observation in \cite{MiPeRe14RGFL} is that $\calL$ defines a GS
$(\PROB(S), \calE,\calR)$ via the explicit representation 
\begin{equation}
  \label{eq:I.RateFcn}
\calL(\rho,\dot \rho)= \calR(\rho,\dot\rho) +
\calR^*(\rho,{-}\rmD\calE(\rho)) + \rmD\calE(\rho)[\dot\rho],
\end{equation}
whenever $\bbQ$ has a unique steady state $\pi\in \PROB(S)$ and $\bbQ$
satisfies the detailed-balance condition with respect to $\pi$ (i.e.\
the Markov process is reversible). We refer to Section \ref{ss:Markov}
for details, where we also highlight that the arising
gradient systems are classical only in the case of diffusion
processes. In case of jumps, one obtains generalized GS involving the
function $\CCC$. In particular, for $\dot p=1{-}2p$ one finds
$([0,1],\calE_\text{Mv},\calR_\text{Mv})$ with $a=\frac12$.

We consider the above stochastic approach as a first microscopic origin
of GS. The second origin involves the concept of evolutionary
$\Gamma$-convergence for GS, see  the
surveys \cite{Serf11GCGF,Miel14?EGCG} for the general ideas.
Here we concentrate on convergence results based on the 
EDP, cf.\ \eqref{eq:I.EDE},  which is an ideal tool 
for doing a limit passage for solutions $u_\eps
:[0,T]\to \bfX$ for a family $(\bfX,\calE_\eps,\calR_\eps)$ of GS
depending on a small parameter $\eps$. The aim is then to derive a 
limiting GS $(\bfX,\calE_0,\calR_0)$ such that a limit $u$ of the
solutions $u_\eps$ is indeed a solution for the limiting GS.
Our Definition \ref{def:EDP-conv} introduces
the concept of EDP-convergence: A family of GS
$(\bfX,\calE_\eps,\calR_\eps)$ converges to the GS
$(\bfX,\calE_0,\calR_0)$ in the EDP sense, if the following holds:
\begin{subequations}\label{eq:I.Def-EDP}
\begin{align}
&   \label{eq:I.Def-a}
\left. \ba{@{\!}c@{\!}} u_\eps:[0,T]\to \bfX \text{ is a}\quad\mbox{}\\
  \text{ solution of\/ }(\bfX,\calE_\eps,\calR_\eps),\\
    u_\eps(0)\weak u^0, \text{ and}\\
 \calE_\eps(0,u_\eps(0)) \to \calE_0(0,u^0){<}\infty \ea \right\} 
 \:  \Longrightarrow  \: 
 \left\{ \ba{@{\!}c@{\!}} \exists\, u \text{ sol.\ of\/ }
   (\bfX,\calE_0,\calR_0) \text{ with  } u(0){=}u^0\\ 
 \text{ and a subsequence } \eps_k\to 0: \\
 \forall\, t\in {]0,T]}{:} \ u_{\eps_k}(t)\weak u(t) \text{ and}\\
 \mbox{}\qquad \qquad \  \calE_{\eps_k}(u_{\eps_k}(t)) \to \calE_0(u(t));
 \ea \right.
\\[0.3em]
&\calE_\eps \Gweak \calE_0 \text{ in }\bfX;\label{eq:I.Def-b} \\[0.3em]
\label{eq:I.Def-c}&
\ba{@{}l} \wt u_\eps(\cdot) \overset{*}{\rightharpoonup}  
 \wt u(\cdot) \text{ in } \rmL^\infty([0,T];\bfX) \text{ and}\\
\sup_{\eps\in {]0,1]},\:t\in [0,T]} \calE_\eps(\wt u_\eps(t)) \leq C < \infty
\ea\!\Big\} \  \Longrightarrow \
 \liminf_{\eps\to 0} {\mathscr D}_\eps(\wt u_\eps) \geq
 {\mathscr D}_0(\wt u) . 
\end{align}
\end{subequations}
When asking only for condition \eqref{eq:I.Def-a} we speak of
pE-convergence, see Definition \ref{def:pE-conv}. Note that \eqref{eq:I.Def-c} enforces a liminf
estimate of De Giorgi's dissipation 
functionals ${\mathscr D}_\eps$ along general functions $\wt u_\eps$,
not only along the solutions of the GS
$(\bfX,\calE_\eps,\calR_\eps)$. Having this liminf estimate, it is
easy to pass to the limit in the $\eps$-dependent energy-dissipation
estimate \eqref{eq:I.EDE}, since the initial energy on the right-hand
side is assumed to converge according to \eqref{eq:I.Def-a}. Then,
applying the EDP for the limiting GS $(\bfX,\calE_0,\calR_0)$ we see
that $u$ is a solution.

In fact, many approaches to evolutionary 
$\Gamma$-convergence establish EDP-convergence, but do not explicitly
state condition \eqref{eq:I.Def-c} as a main result. E.g.\ the
Sandier-Serfaty approach \cite{SanSer04GCGF,Serf11GCGF}, where the
terms $\int_0^T\calR_\eps\dd t$ and $\int_0^T\calR_\eps^*\dd t$ are
treated separately, provides EDP-convergence. Our approach is more
general than the latter, since we only ask that the sum
$\int_0^T\calR_\eps\dd t+ \int_0^T\calR_\eps^*\dd t$ behaves well, but
not necessarily the individual terms. This has two effects: (i) we can
allow for general functions $\wt u_\eps$, and (ii) it can lead to
exchanges between the two terms in the limit $\eps\to 0$.  Point (i)
is important to explore ${\mathscr D}_0$ outside of the set of
solutions and thus providing the full information about the GS
$(\bfX,\calE_0,\calR_0)$, while the set of solutions of the limit
equation $\dot u=V_0(t,u):= \rmD_{\dot u}\calR_0 (u,{-}\rmD_u
\calE_0(t,u))$ only contains information on $V_0$.
Point (ii) is relevant for another important message of this
paper. The EDP-limit of classical GS can be a generalized GS. This
phenomenon is considered as another microscopic origin of generalized
GS.  

Here we provide three different examples for point (ii), the
first of which is an ODE example in Section \ref{sss:EE+quad}, while
Sections \ref{se:Membrane} and \ref{se:D2R} contain more elaborate
examples treating the membrane limit of a thin-layer and the limit
of diffusion to reaction, respectively.

For the membrane limit we consider a diffusion equation with a thin
layer with very small diffusivity. In \cite{Lier12VME,Lier13PBBS}
pE-convergence to the membrane limit was established; however 
EDP-convergence was not studied. We start with the diffusion equation 
which is the gradient flow for the classical GS $(\PROB(\Omega), 
\calE, \calR_\eps^*)$ with $\calE(u)=\int_\Omega u\log(2 u) \dd x$
and 
$\calR_\eps^*(u,\xi) = \frac12 \int_\Omega a_\eps(x) (\pl_x\xi)^2 u \dd x $.
Using suitable scalings for the diffusion coefficient $a_\eps(x)$
Theorem \ref{th:Membrane} provides  EDP-convergence to the
generalized GS $(\PROB(\Omega), \calE, \calR_0^*)$ with 
\[
\calR_0^*(u,\xi)= \int_{-1}^0 \frac{a}2 (\pl_x\xi)^2 u\dd x
  + a_* \sqrt{u(0^-)u(0^+)}\: \CCC^*( \xi(0^+){-}\xi(0^-))
+ \int_0^1\frac{a}2 (\pl_x\xi)^2u \dd x ,
\]
where $u(0^-)$ and $u(0^+)$ denote the limit of $u(x)$ at $x=0$ 
from the left and from the right, respectively. Thus, 
$\calR_0^* $ involves $\CCC^*$ and is therefore non-quadratic. 

Section \ref{se:D2R} follows
\cite{PeSaVe10DRGC,PeSaVe12CRGL,AMPSV12PLWG} by considering the limit
from pure diffusion in physical space and along a reaction-path
variable $y\in \Upsilon\subset \R$ to a limit of a reaction-diffusion
system on $\Omega$. The Fokker-Planck equation reads
\[
\dot u = m_\Omega \Delta_x u + \tau_\eps \pl_y \Big( \pl_y u +
\frac1\eps u\,\partial_yV(y)  \Big), 
\]
where $V$ is a potential with two global minima $y_0$ and $y_1$ 
and one global maximum in-between.  
This equation is generated by the classical GS
$(\PROB(\Omega{\ti}\Upsilon), \calE_\eps,\calR^*_\eps)$ where $\calE_\eps$ is the relative entropy and $\calR_\eps^*$ is the quadratic Wasserstein dissipation potential, see \eqref{eq:D2R.F.R}.
Theorem \ref{th:D2R-calJ} establishes EDP-convergence to a generalized
GS $(\PROB(\Omega\ti\{y_0,y_1\}),\bfE,\bfR^*)$, where $\bfR^*$ again involves 
the non-quadratic function $\CCC^*$.

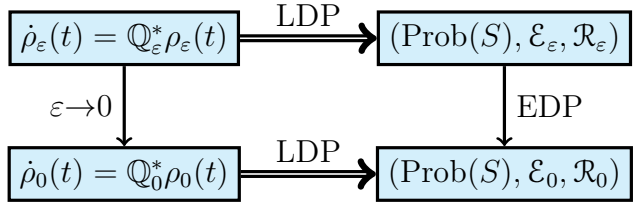
\begin{figure}
\centering 
\begin{tikzpicture}[scale=1.0]
\tikzstyle{every node} = [draw, very thick, fill=cyan!15];
\node (KFeps) at (0,1.8) {$\dot \rho_\eps(t)= \bbQ^*_\eps \rho_\eps(t)$};
\node (GSeps) at (5,1.8) {$(\PROB(S),\calE_\eps, \calR_\eps)$};
\node (KF0) at (0,0) {$\dot \rho_0(t)= \bbQ^*_0 \rho_0(t)$};
\node (GS0) at (5,0) {$(\PROB(S),\calE_0, \calR_0)$};
\tikzstyle{every node} = [];
\draw[->,double,very thick,very thick] (KFeps)--(GSeps) node[pos=0.5,above]{LDP};
\draw[->,double,very thick] (KF0)--(GS0) node[pos=0.5,above]{LDP};
\draw[->,very thick] (KFeps) -- (KF0)
node[pos=0.5,left]{$\eps{\to} 0$};
\draw[->,very thick] (GSeps) -- (GS0) node[pos=0.5,right]{EDP};
\end{tikzpicture}
\caption{For reversible, time-continuous Markov processes the
  large-deviation principle (LPD) of Section
  \ref{ss:Markov} provides a (generalized) gradient structure. This
  mapping commutes with taking the limit $\eps\to 0$ and
  EDP-convergence, respectively.}
\label{fig:LDPvsEDP}
\end{figure}

We conclude our introduction by a general and surprising
observation. The three main models in this work (i.e.\ the ODE,
the membrane, and the reaction-to-diffusion model in Sections
\ref{sss:EE+exp}, \ref{se:Membrane}, and \ref{se:D2R}, respectively)
can be seen as Kolmogorov forward equations for  naturally
associated Markov processes. Thus, the large-deviation theory of
Section \ref{ss:Markov} is applicable and provides entropic GS
$(\PROB(S),\calE_\eps, \calR_\eps)$ for the associated Kolmogorov
forward equations $\dot\rho_\eps = \bbQ^*_\eps\rho_\eps$ for each
$\eps\in {]0,1[}$ as well as for $\eps=0$. The limit for $\eps
=0$ can be also defined in terms of the classical convergence for
Markov processes asking $\rho_\eps(t)=\ee^{t\bbQ^*_\eps} \rho(0) \weaks
\rho(t)=\ee^{t\bbQ^*_0} \rho(t)$. 
Ignoring the (linear) Markovian structure, we can also consider 
EDP-convergence of the induced entropic GS. In all our three examples
we find the surprising result that
the EDP-limit is exactly the entropic GS of the limiting Markov
process. This means that applying the described large-deviation
principle and taking the limit $\eps\to 0$ (either on the level of 
Markov semigroups or as EDP-convergence for GS) commute, see Figure
\ref{fig:LDPvsEDP}. 
This result appears naturally, if we use representation
\eqref{eq:I.RateFcn} of the rate function $\mathscr I$ giving
\[
\mathscr I(\rho)=\int_0^T \calR(\rho,\dot
\rho)+\calR^*(\rho,{-}\rmD\calE(\rho)) +\rmD\calE(\rho)[\dot\rho] \dd
t = {\mathscr D}(\rho)+ \calE(\rho(T))-\calE(\rho(0)).
\]
Hence, the above large-deviation principle exactly encodes the
energy-dissipation principle, and EDP-convergence for the induced
entropic GS can be interpreted as $\Gamma$-convergence of the rate
functionals.

The question how general this observation about the interchangeability
of the suitable large-deviation principles and the EDP-convergence is,
seems to be challenging, but goes beyond the scope of this work. We
mention that in \cite{BonPel14?QRIL} similar relations between
large-deviation principles and evolutionary $\Gamma$-convergence are
studied. 

As a final general remark, we emphasize that this paper focuses
on the modeling aspects of the emergence of
generalized GS. Thus, we do not give the full analytical details in
terms of estimates and convergences in the proper functional
spaces, but rather highlight the structures and manipulations needed to
understand the corresponding limit procedures.

\section{Classical and generalized gradient systems}
\label{se:gGS}

We now convert the formal ideas from the introduction into rigorous
mathematical statements. We call a triple $(\bfX,\calE,\calR)$ a
gradient system (GS), if $\bfX$ is a Banach
space, $\calE:[0,T]\ti \bfX\to \R_\infty:=\R{\cup} \{\infty\}$ is a
functional (such as the free energy, the negentropy, etc.), and
$\calR:\bfX\ti \bfX\to [0,\infty]$ is a 
dissipation potential, which means that for all $q\in \bfX$
the functional $\calR(u,\cdot):\bfX\to \R_\infty$ is lower
semicontinuous, nonnegative, convex, and satisfies $\calR(u,0)=0$. 
In this section, we allow for the case that the energy functional
depends on the time variable $t\in [0,T]$ to show that the abstract
principle is valid in this general case. However, for notational
convenience we will restrict to the autonomous case (i.e. 
$\pl_t\calE(t,u)\equiv 0$) in all other parts.  

We speak of a \emph{classical GS}, if
$\calR(u,\cdot)$ is quadratic, i.e.\ there exists a symmetric and
positive definite operator $\bbG$ such that $\calR(u,v)=
\frac12\langle \bbG(u)v,v\rangle$. 
However, plasticity requires
non-quadratic dissipation potentials, e.g.\ of the form
$\calR(\dot\pi)= \sigma_\text{yield}\|\dot\pi\|_{\rmL^1} +
\tfrac12\mu_\text{visc}\|\dot\pi\|_{\rmL^2}^2$, see
\cite{Miel03EFME,MieRou08?RIST}.  In particular, the
rate-independent case is based on $\calR(u,\lambda v)=\lambda \calR(u,v)$
for all $\lambda>0$, which is incompatible with a quadratic form. If
$\calR(u,\cdot)$ is non-quadratic, we call $(\bfX,\calE,\calR)$ a
\emph{generalized GS}.

\subsection{Variational principles for gradient systems}
\label{ss:VarPrinc}

The following proposition from convex analysis shows that there are
several completely equivalent formulations of the generalized force
balance \eqref{eq:I.FB}. The equivalences of the points (ii) to (iv)
below are also called Fenchel equivalences, cf.\ \cite{Fenc49CCF}. The
essential tool is the Legendre-Fenchel transform $\Psi^*:\bfX^*\to
\R_\infty$ of a convex function $\Psi:\bfX\to \R_\infty$ defined via
\[
\Psi^*(\xi) := \sup\set{\langle \xi,v\rangle-\Psi(v) }{ v\in \bfX}. 
\] 
In a reflexive Banach space we have $(\Psi^*)^*=\Psi$. 

\begin{proposition}[Equivalent formulations]\label{pr:Fenchel}
  Let $\bfX$ be a reflexive Banach space and $\Psi:\bfX\to \R_\infty$ be
  proper, convex, and lower semicontinuous. Then, for every $\xi \in
  \bfX^*$ and every $v\in \bfX$ the following five
  statements are equivalent:
\begin{align*}
&\text{\ (i) \ }\ v \in \mathop{\mathrm{Arg\,min}}\limits_{w\in \bfX}
             \big(\Psi(w) - \langle\xi,w\rangle\big); \qquad 
\text{(ii) } \xi \in \pl\Psi(v); \\ 
&\text{(iii) } \ \Psi(v) +\Psi^*(\xi) = \langle  \xi,v\rangle; \\[0.5em]
&\text{(iv) } \ v \in \pl\Psi^*(\xi); \qquad  \qquad 
\text{(v) } \ \xi \in \mathop{\mathrm{Arg\,min}}\limits_{\eta\in \bfX^*}
\big(\Psi^*(\eta)-\langle \eta,v\rangle\big). 
\end{align*}
\end{proposition}
Note that the definition of $\Psi^*$ immediately implies the
Young-Fenchel inequality $\Psi(w)+\Psi^*(\eta)\geq \langle 
\eta,w\rangle$ for all $w$ and $\eta$. Thus, (iii) expresses an
optimality as well. 

Defining the dual dissipation potential $\calR^*$ via
$\calR^*(u,\cdot):= (\calR(u,\cdot))^*$ we can apply these equivalences
to reformulate \eqref{eq:I.FB} in the following ways:

\begin{description}\itemsep0.5em
\item[\ (I) \  
Rayleigh principle] \mbox{}\hfill
  \cite{
Rayl1871GTRV} 
\\[0.3em]
\hspace*{2em}
(RP)\quad $\dot u \in \mathop{\mathrm{Arg\,min}}\limits_{v\in\bfX}\Big( \calR(u,v) 
- \langle \rmD  \calE(t,u), v\rangle \Big)$;

\item[(II) \  Force balance in $\bfX^*$]\mbox{}\hfill Rayleigh-Biot
  equation \cite{Rayl1871GTRV,Biot55VPIT} 
\\[0.3em]
\hspace*{2em}
(FB)\ \quad $ 0 \in \pl_{\dot u} \calR(u,\dot u) + \rmD\calE(t,u) \ \in \bfX^*$; 

\item[(III) Power balance in $\R$] \mbox{}\hfill De Giorgi's $(\calR,\calR^*)$
  formulation \cite{DeMaTo80PEMS} 
\\[0.3em]
\hspace*{2em}
(PB)\ \quad  $  \calR(u,\dot u) + \calR^*(u,{-}\rmD
  \calE(t,u) )=  -\langle\rmD\calE(t,u),\dot u \rangle$;

\item[(IV) Rate equation in $\bfX$] \mbox{}\hfill Onsager equation
  \cite{Onsa31RRIP}
  \\[0.3em]
  \hspace*{2em} (RE)\ \quad $\dot u \in \pl_\xi \calR^*(u,{-}
  \rmD\calE(t,u) )\ \in \bfX$;

\item[(V) \ Maximum dissipation principle] \mbox{} \hfill cf.\ e.g.\
  \cite{HacFis08RBPM} 
\\[0.3em]
\hspace*{2em}
(MDP)\quad  $ \rmD\calE(t,u) \in \mathop{\mathrm{Arg\,max}}\limits_{\xi\in\bfX^*}\Big( \langle
\xi,\dot u\rangle - \calR^*(u,\xi)\Big)$.
\end{description}
In fact, \cite[Eqn.\,(26)]{Rayl1871GTRV} also includes the
kinetic energy $\calT$, which we omit in our approximation, namely 
$\frac\rmd{\rmd t} \big( \rmD_{\dot u}\calT(u,\dot u)\big) + \rmD_{\dot u}
\calR(u,\dot u) + \rmD_q\calE(t,u)=0$.

Before returning to the general situation, we highlight the three different
cases (II)--(IV) for the classical viscous dissipation, i.e.\
$\calR(u,v)=\frac12\langle \bbG(u) v,v\rangle$ and
$\calR^*(u,\xi)=\frac12\langle \xi,\bbK(u)\xi\rangle$ with
$\bbK(u)=\bbG(u)^{-1}$. Then, we have 
\begin{align*}
\text{(FB) }\quad  &\bbG(u) \dot u = - \rmD \calE(u)\qquad \quad
\text{(RE) }\quad \dot u = - \bbK(u)\rmD \calE(u)=: - \nabla_{\!\bbG}\calE(u)\\[0.5em]
\text{(PB) }\quad &\frac12\langle \bbG(u) \dot u, \dot u\rangle + 
\frac12\big\langle \rmD \calE(u),\bbK(u)\rmD \calE(u)\big\rangle= - \langle \rmD
\calE(u),\dot u\rangle,
\end{align*}
where (RE) can be seen as a ``\emph{gradient-flow equation}'', as
$\nabla_\bbG$ is the gradient operator. 

\subsection{The energy-dissipation  principle}
\label{ss:EDP}

The above formulations can already be understood in a variational sense,
since the evolution is expressed by extremizing a functional or by 
variations or derivatives of the two functionals $\calE$ and
$\calR$. However, for mathematical purposes it is desirable to have
formulations in terms of a minimization problem 
for the whole solution trajectories
$u:[0,T]\to \bfX$. One such principle can be derived on the basis of
the power balance (PB) by integration in time and using the chain rule
and finally employing the 
Young-Fenchel inequality $\Psi(w)+\Psi^*(\eta)\geq \langle
\eta,w\rangle$, cf.\ \cite{DeMaTo80PEMS} or the survey
\cite{Miel14?EGCG}. This leads to the celebrated energy-dissipation
principle, also called De Giorgi's $(\calR,\calR^*)$ principle, see
\cite{AmGiSa05GFMS} for some historical remarks. 

\begin{theorem}[De Giorgi's energy-dissipation principle]\label{th:EDP}
Under suitable technical conditions on $(\bfX,\calE,\calR)$ a function
$u:[0,T]\to \bfX$ satisfies (I)--(V) from above for almost all $t\in [0,T]$
if and only if the \emph{energy-dissipation balance (EDB)} holds: 
\[
\mafo{(EDB)} \qquad \left\{ 
\begin{aligned} 
\calE(T,u(T)) +
  \int_0^T \calR(u,\dot u)& + \calR^*(u,{-}\rmD\calE(t,u) ) \dd t\\[-0.6em]
&= \ \calE(0,u(0)) + \int_0^T \pl_t \calE(t, u(t)) \dd t. 
\end{aligned} \right.
\]
Under additional technical conditions it is sufficient to have
only 
the upper estimate where ``$=$'' is replaced by ``$\leq$''. In this
case, we speak of the \emph{energy-dissipation estimate (EDE)}.
\end{theorem}


\subsection{Examples of generalized gradient structures}
\label{ss:ExamGGS}

Here we give some examples of generalized gradient structures. First,
we discuss dissipative material models like plasticity or shape-memory
materials that  form a huge class of generalized GS. Second, we 
treat nonlinear reaction-diffusion systems (RDS), which will be closer
to the main theme of this paper. The third class of examples concerns
reversible Markov processes, where the Kolmogorov
forward equation has a gradient structure with the relative entropy as
energy functional. This latter class is so important that it is
treated in the subsequent Subsection \ref{ss:Markov}. 

\subsubsection{Dissipative material models} 
\label{sss:DissMater}

The state of a body $\Omega\subset \R^d$, composed of so-called 
dissipative materials (also called standard generalized
materials), 
is given in terms of the elastic
displacement $\bfu:\Omega\to \R^d$ and an additional internal variable
$z:\Omega\to \R^k$. The latter may describe plastic deformations,
damage, phase-field variables, magnetization, or other internal states of the
material. The total stored energy $\calE$ depends on $\bfu$, $z$, and
usually also on a process time $t\in[0,T]$, i.e.\ $\calE(t,\bfu,z)$. 
As introduced in the theory of standard generalized materials in
\cite{HalNgu75MSG} the dissipative forces are given in terms of a
(primal) dissipation potential $\calR$ that also may include
viscoelastic terms:
\[
\calR(\dot\bfu,\dot z)= \calR_\mafo{diss}(\dot z) +
\calR_\mafo{visc}(e(\dot\bfu)) ,
\]
where $e(\dot\bfu)=\frac12(\nabla \dot \bfu
{+}(\nabla\dot\bfu)^\top)$. As before, the corresponding force
balance equations (FB) are  
\[
0= \rmD_{\dot\bfu}\calR(\dot\bfu,\dot z) + \rmD_\bfu \calE(t,\bfu,z),
\quad 
0\in\pl_{\dot z} \calR(\dot\bfu,\dot z) + \rmD_z \calE(t,\bfu,z)
\]
with $\pl_{\dot z}\calR(\dot\bfu,\dot z)$ denoting the set-valued convex subdifferential.

While the viscoelastic potential
$\calR_\mafo{visc}$ is assumed to be quadratic in many applications,
the potential $\calR_\mafo{diss}$ for the internal variables $z$ is
often supposed to be non-quadratic. E.g.\ in viscoplasticity with
yields stress $\sigma_\mafo{yield}$ one takes the form 
\[
\calR_\mafo{diss}(\dot z)=\int_\Omega \Big(\sigma_\mafo{yield} |\dot z|_1+
\nu |\dot z|^{1+\delta} \Big) \dd x,
\]
where $\delta>0$ is usually taken small, e.g.\ in $\delta=0.012$ in
\cite{ZRSRZ96TDIT}. The weak growth of order $1{+}\delta$ is sometimes  even
replaced by a growth $O(|\dot z|\log|\dot z|)$ as given by our
function $\CCC$ (see
e.g.\ \cite[Eqn.\,(5)]{RoRaGo00WHHA} and \cite{BonPel14?QRIL}). 

For later reference, we mention the very simple scalar hysteresis
model of a so-called play operator. It is given by the generalized GS
$(\R,\calE_\mafo{play},\calR_\mafo{play})$ with 
\begin{equation}
  \label{eq:Play}
  \calE_\mafo{play}(t,z)=\frac12 z^2 -\ell(t) z \quad \text{and} \quad 
\calR_\mafo{play}(\dot z)= r|\dot z| \text{ with }r>0.
\end{equation}
It serves as a limit for evolutionary $\Gamma$-convergence in Example
\ref{ex:pEnotEDP} as well as a
large-deviation limit in \cite{BonPel14?QRIL}.

\subsubsection{Nonlinear reaction-diffusion systems} 
\label{sss:nonlRDS}

We consider concentrations $\bfc(t):\Omega \to
{[0,\infty[}^I$ of chemical species $C_1,\ldots C_I$ that can react according
to $R$ reactions of mass action type given by a stoichiometric
relation 
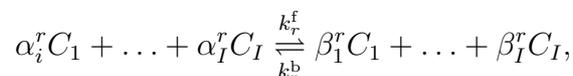
\[
\alpha_i^r C_1 + \ldots +\alpha^r_I C_I
\overset{k_r^\rmf}{\underset{k^\rmb_r}{\rightleftharpoons}} 
\beta^r_1 C_1+ \ldots +\beta^r_I C_I,
\]
where $r=1,\ldots, R$ is the index of the reaction, $k^\rmf_r$ and
$k^\rmb_r$ are the forward and backward reaction coefficients, and the
stoichiometric coefficients $\alpha^r_i$ and $\beta^r_i$ are
nonnegative integers.  The reaction-diffusion system (RDS) for the
concentrations $\bfc=(c_1,\ldots,c_I)$ takes the form
\begin{equation}
  \label{eq:nRDS1}
  \dot \bfc = \bbD \Delta \bfc - \bfR(\bfc) \  \text{ with } 
\bfR(\bfc) := \sum_{r=1}^R \big(k_r^\rmf
  \bfc^{\bfalpha^r} {-} k^\rmb_r\bfc^{\bfbeta^r} \big)\big(\bfalpha^r{-}
  \bfbeta^r\big)
\end{equation}
and $\bbD=\diag(\delta_i)_{i=1,\ldots,I}$, where $\delta_i>0$. With the stoichiometric
vectors $\bfalpha^r = ( \alpha_i^r)_i$ and $\bfbeta^r=(\beta^r_i)_i
\in \bbN_0^I$ we define the monomials in the form $\bfc^\bfalpha:=
\prod_{i=1}^I c_i^{\alpha_i}$.

It was shown in \cite{Miel11GSRD,Miel13TMER} that \eqref{eq:nRDS1} has
a (classical) gradient structure under the additional assumption of
the detailed-balance condition, which means that
\begin{equation}
  \label{eq:nRDS.DBC}
  \exists \,\bfw=(w_i)_i: \quad w_i>0 \ \text{ and } \ k^\rmf_r
  \bfw^{\bfalpha^r} = k^\rmb_r \bfw^{\bfbeta^r}\text{ for }r=1,\ldots,R.
\end{equation}
Using the Boltzmann function $\lambda_\rmB(z)=z \log z -z +1$ 
we define the relative entropy 
\[
\calE(\bfc)=\int_\Omega \sum_{i=1}^I \lambda_\rmB(c_i(x)/w_i) w_i \dd
x,
\]
which gives rise to the vector of thermodynamic driving forces (also
called chemical potentials) $\bfmu=\rmD \calE(\bfc)$ with $\mu_i=
\log(c_i/w_i)$. Because of the logarithm laws they satisfy the relation 
$\bfalpha^r \cdot \rmD\calE(\bfc)= \sum_{i=1}^I \alpha^r_i
\log(c_i/w_i) = \log(\bfc^{\bfalpha^r}) - \log(\bfw^{\bfalpha^r})$.
Thus, using the detailed-balance conditions we obtain the relation 
\begin{equation}
  \label{eq:nRDS3}
  \big(\bfalpha^r{-} \bfbeta^r\big) \cdot ({-}\rmD\calE(\bfc)) =
 \log(k^\rmb_r \bfc^{\bfbeta^r}) -\log(k^\rmf_r\bfc^{\bfalpha^r} ).
\end{equation}
To construct the dual dissipation potentials we may choose any scalar,
strictly convex dual dissipation functional $\psi:\R\to \R$ with
$\psi(0)=\psi'(0)=0$ and $\psi''(0)>0$ and let
\begin{align*}
&\calR^*(\bfc, \bfmu):= \int_\Omega \Big[\frac12 \sum_{i=1}^I \delta_i c_i
|\nabla \mu_i|^2 \ + \ \sum_{r=1}^R H^r_\psi(\bfc)
\psi\big((\bfalpha^r{-} \bfbeta^r) {\cdot} \bfmu \big)\Big] \dd x \\
&\text{with } H^r_\psi(\bfc):=   \frac{
  k^\rmb_r \bfc^{\bfbeta^r}\;-\;\: k^\rmf_r \bfc^{\bfalpha^r} }{\psi'\big( \log (k^\rmb_r
  \bfc^{\bfbeta^r}) {-}\log( k^\rmf_r \bfc^{\bfalpha^r}) \big)}.
\end{align*}

Using $\rmD_\bfmu\calR(\bfc,\bfmu)=
  -\big(\Div (\delta_ic_i \nabla \mu_i)\big)_i + \sum_{r=1}^R
H^r_\psi(\bfc)\psi'\big((\bfalpha^r{-} \bfbeta^r ) {\cdot}
\bfmu\big)\big(\bfalpha^r{-}\bfbeta^r\big)$ and \eqref{eq:nRDS3}, we
easily see that the nonlinear RDS \eqref{eq:nRDS1} satisfying the
detailed-balance condition is 
generated by the pair $\calE$ and $\calR$, i.e.\ $
\dot \bfc= \bbD\Delta \bfc -\bfR(\bfc) = \rmD_\bfmu \calR \big(\bfc,
{-}\rmD\calE(\bfc)\big)$.  
Thus, we have found a family of generalized gradient structures of
the nonlinear RDS \eqref{eq:nRDS1}. 

The case of quadratic
$\calR(\bfc;\cdot ) $, i.e., $\psi(\eta)=\eta^2/2$ was introduced in
\cite{Miel11GSRD}, leading to the logarithmic means 
$H^r_\text{quadr}(\bfc)= \Lambda( k^\rmf_r \bfc^{\bfalpha^r},k^\rmb_r
\bfc^{\bfbeta^r})$ with $\Lambda(a,b)=(a{-}b)/(\log a{-}\log b)$, see
also \cite{Maas11GFEF,ErbMaa12RCFM,MaaMie15?GSRC1}. However, it was
already criticized in the 1930s that the linear relation $\dot\bfc = -
\bbK(\bfc)\bfmu$ (i.e.\ $\calR^*(\bfc,\bfmu)=\frac12\langle
\bfmu,\bbK(\bfc)\bfmu\rangle$ is quadratic) arising from Onsager's
principle is not suitable for 
chemical reactions if one wants to model systems that are not very
close to thermal equilibrium. As a solution Marcelin and de Donder
introduced exponential dependencies between $\bfmu$ and $\dot\bfc$, 
see \cite[Def.\,3.3]{Fein72CKCC} or \cite[Eqn.\,(11)]{GKZD00RDRK}. In  
\cite{Grme10MENT} Remark iii on p.\,77 gives some historical comments
and Eqn.\,(69) explicitly features an exponential dissipation potential $\Xi$
involving the function $\big(\ee^{\xi/2} + \ee^{-\xi/2}-2)$. Since the
choice $\psi(\xi)=\CCC^*(\xi)$ is
central for our paper, we give $\calR^*$ explicitly for this case, viz.\
\begin{equation}
  \label{eq:nRDS-R*CCC}
  \calR^*(\bfc,\bfmu) =\int_\Omega \Big[
\sum_{i=1}^I \frac{\delta_i c_i}2 |\nabla c_i|^2 \, + \sum_{r=1}^R 
\sqrt{k^\rmf_r\bfc^{\bfalpha^r} k^\rmb_r\bfc^{\bfbeta^r}} \,
\CCC^*\big((\bfalpha^r{-}\bfbeta^r)\cdot \bfmu  \big)   \Big] \dd x.   
\end{equation}
We will see that exactly the same structure, up to a trivial scaling
factor $1/2$, arises via the large-deviation principle described next,
see also \cite{MPPR15?LDPC}.

\subsection{Markov processes, large deviations, and GS}
\label{ss:Markov}

Here we give a rough sketch of the theory in \cite{MiPeRe14RGFL} about
 gradient structures for the Kolmogorov forward equation $\dot
\rho=\bbQ^* \rho$ of Markov
processes satisfying a detailed-balance condition, which are also
called reversible Markov processes, for short. The idea that
large-deviation principles generate gradient structures goes back to
\cite{OnsMac53FIP} (see Eqn.\,(4-21) therein for a quadratic version of the
energy-dissipation principle derived by large-deviations, called
Boltzmann's principle).  The mathematical theory was developed
only recently, see \cite{ADPZ11LDPW,ADPZ13LDGF,MiPeRe14RGFL}.  

In Section \ref{sss:RevMP-GS} we first describe a time-dependent
large-deviation principle for general Markov processes 
providing a formula
for the rate function $\mathscr I(\rho(\cdot)) =\int_0^T\calL(\rho(t),\dot\rho(t))
\dd t$ and then present the result of
\cite{MiPeRe14RGFL}, which shows that for reversible Markov processes
the functional $\calL$ is induced by an EDP for a GS
$(\PROB(S),\calE,\calR)$. In Sections \ref{sss:FiniteMP} to
\ref{sss:MP-Memb} we then discuss a few applications of the abstract
result in Theorem \ref{th:LDgivesGS}.  

\subsubsection{Gradient structures obtained via large deviations}
\label{sss:RevMP-GS}
We consider a compact metric space $S$ and denote by $\PROB(S)$ the
subset of probability Radon measures on $S$ equipped with the narrow
convergence $\weaks$ defined by duality with continuous, bounded
functions. The Kolmogorov forward equation
$\dot \rho = \bbQ^* \rho$ describes the evolution of the probability
laws $\rho(t)$ of a Markov process $(X_t)_{t\geq 0}$, 
if the law of $X_0$ is given by $\rho(0)\in
\PROB(S)$. The Markov generator is given as $\bbQ$ acting on
functions on $S$, while its dual $\bbQ^*$ acts on measures such that 
$\int_S (\bbQ f) \dd \rho= \int_S f \dd (\bbQ^*\rho)$. 

Considering $N$ independent realizations $(X^{(i)}_t)_{t\geq 0} $,
$i=1,\ldots,N$ of the underlying Markov process, the
measure-valued empirical process $\rho^N(t):=\frac1N \sum_{i=1}^N
\delta_{X_t^{(i)}} \in \PROB(S)$ can be defined. Using the law of large numbers the
limit $N\to\infty$ gives $\rho^N(t)\weaks \rho(t)$, which solves the
Kolmogorov forward equation $\dot \rho= \bbQ^*\rho$, see e.g.\
\cite[Thm.\,2.3.1]{Reng13MIWG}. Under suitable assumptions, see
\cite{MiPeRe14RGFL}, it is shown in \cite{FenKur06LDSP} that the
empirical process $\rho^N$ satisfies a large-deviation principle with
a rate function $\mathscr I(\rho(\cdot))$, i.e.\
\[
\bbP\big( \rho^N(\cdot) \approx \wh\rho(\cdot)\big)\  \simeq \ 
\ee^{-N \mathscr I(\wh \rho(\cdot))}, 
\]
see the above references for the proper definition of ``$\simeq$''.
The main result is that $\mathscr I$ has the form $ \mathscr I(\rho(\cdot))=
\int_0^T \calL(\rho(t),\dot \rho(t)) \dd t$, where $\calL$ is the
Legendre transform of the functional $\calH(\rho,\cdot)$ (\,i.e.\
$\calL(\rho,v)=\sup \int_\Omega \xi \rmd v - \calH(\rho,\xi)\:$) given
via the explicit formula
\[
\calH(\rho,\xi) := \int_S \ee^{-\xi(s)}(\bbQ \,\ee^{\xi})(s) \dd
\rho(s). 
\]
We emphasize the simplicity of this formula and the (separate)
linearity in $\rho$ and in $\bbQ$.

A main observation in \cite{MiPeRe14RGFL} is that the deterministic
case, which is given by the relation
\[
\mathscr I(\rho(\cdot))= 
\int_0^T \calL(\rho(t),\dot\rho(t))\dd t=0,
\]
can be interpreted as an energy-dissipation principle if and only if
the Markov process is reversible, which is the same as asking for the
detailed balance condition (cf.\ \eqref{eq:nRDS.DBC}) for the linear
Kolmogorov forward equation $\dot \rho = \bbQ^*\rho$. 
Hence, we now further assume that there
exists a stationary measure $\pi\in \PROB(S)$ which has, without loss of
generality, the full set $S$ as its support. We say that \emph{$\bbQ$
satisfies the detailed balance condition with respect to $\pi$}, if
\begin{equation}
  \label{eq:DBC-Markov}
 \int_S f \,\big(\bbQ\,g\big) \dd \pi = \int_S g\,\big(\bbQ\,f) \dd \pi  
\end{equation}
for all $f$ and $g$ in the domain if $\bbQ$. Choosing $g\equiv 1$, we
find that $\bbQ^*\pi=0$, i.e.\ the detailed balance condition implies
the stationarity. 

This version of the detailed-balance condition for Markov processes 
coincides with the detailed-balance condition for chemical reactions
in \eqref{eq:nRDS.DBC}. Indeed, if the ODE case $\dot \bfc=-\bfR(\bfc)$ is
 linear, i.e.\ $\dot\bfc =\bbA\bfc$ with $\bbA \bfw=0$, then
\eqref{eq:nRDS.DBC} means $\bbA_{ij}w_j = \bbA_{ji} w_i$. Setting
$S=\{1,\ldots,I\}$ and $\bbQ=\bbA^*$ gives \eqref{eq:DBC-Markov}.

In the sequel we will use the Radon-Nikodym derivative of $\rho$ with
respect to $\pi$ denoted by $f = \frac{\rmd\rho}{\rmd\pi}\in
\rmL^1_{\geq 0}(S,\pi)$ and defined via \ $\int_B 1 \dd\rho =
\int_B f\dd \pi$ for all Borel sets $B\subset S$.

\begin{theorem}[{\cite[Sec.\,3]{MiPeRe14RGFL}}]\label{th:LDgivesGS}
If the Markov process $(X_t)_{t\geq0}$ on $S$ is reversible, i.e.\ a
stationary measure $\pi\in \PROB(S)$ satisfying the detailed-balance condition
\eqref{eq:DBC-Markov} exists for the Kolmogorov forward equation
$\dot \rho = \bbQ^*\rho$, then the large-deviation rate functional
$\int_0^T\calL(\rho,\dot\rho) \dd t$ has the form of an
energy-dissipation principle, namely
\[
\calL(\rho,v)= \calR(\rho,v) + \calR^*(\rho,-\rmD\calE(\rho)) +\int_S
\rmD\calE(\rho) \dd v,
\] 
where the gradient structure $(\PROB(S),\calE,\calR^*)$ is given by $
\calE(\rho)= \frac12\int_S \lambda_B\big( \frac{\rmd \rho}{\rmd
  \pi}\big) \dd \pi $ and the dual dissipation potential
\begin{equation}
  \label{eq:calR-LDP}
\calR^*(\rho, \xi)= \int_S \Big((\sqrt{f}\:\ee^{-\xi}) \big(\bbQ(\sqrt{f}
\:\ee^{\xi})\big) - \sqrt f \big(\bbQ \sqrt f\big) \Big) \dd \pi,
\quad \text{where }f = \frac{\rmd\rho}{\rmd\pi}.
\end{equation}
\end{theorem}

The cited reference contains not only a full proof, but also specifies
under what assumptions this implication is in fact an equivalence,
i.e.\ the existence of a gradient structure implies the existence of a
steady state satisfying the detailed-balance condition.

Since the arguments and proofs in \cite{MiPeRe14RGFL} are quite
involved, we highlight here the main structures and formal
calculations to see that $(\PROB(S),\calE,\calR^*)$ is a GS 
and that it generates the Kolmogorov equation $\dot
\rho=\bbQ^*\rho$. 

We first observe that  $\calR^*$ is defined in terms of $\calH$ via 
\[
\calR^*(\rho,\xi)= \calH\Big(\rho,\xi{+}\frac12\log f\Big) -
\calH\Big(\rho,\frac12 \log f\Big), \quad \text{where }\frac12 \log f =
\rmD\calE(\rho) = \frac12\,\log\frac{\rmd\rho}{\rmd\pi}. 
\]
Obviously, $\calR^*(\rho,\cdot)$ is convex if and only if
$\calH(\rho,\cdot)$ is convex. The latter is independent 
of the detailed-balance condition and can be established as
follows.
Consider the Markov semigroup $P_t=\ee^{t\bbQ}$ for $t\geq 0$.  For
fixed $\rho\in \PROB(S)$ and $t\geq 0$ define
\[
\calA_t(\xi):=\int_S \ee^{-\xi}\,P_t(\ee^\xi) \rmd \rho = \int_{S\ti S}
\ee^{-\xi(x)+ \xi(y)} p_t(x,\rmd y)\rho(\rmd x) ,
\]
where $p_t$ with $p_t(x,\cdot) \in \PROB(S)$ denotes the time-dependent
Markov kernel. From the convexity of $\xi\mapsto \ee^{-\xi(x)+
  \xi(y)}$ and the nonnegativity of $p_t $ and $\rho$, we conclude
that $\xi \mapsto \calA_t(\xi)$ is convex. Using $\frac1t(P_t \eta-\eta)\to 
\bbQ \eta$ and $\calA_0(\xi)\equiv 1$, we see that 
\[
\calH(\rho,\xi)\ = \ \lim_{t\to 0} \frac1t(\calA_t(\xi)-\calA_0(\xi))
\ = \ \lim_{t\to 0} \frac1t(\calA_t(\xi)-1)
 \]
is also convex in $\xi$.

By definition we have  $\calR^*(\rho,0)=0$, and
the detailed-balance condition implies the time reversibility
$\calR^*(\rho,{-}\xi)=\calR^*(\rho,\xi)$. Since this implies
$\rmD\calR^*(\rho,0)=0$, convexity gives the positivity 
$\calR^*(\rho,\xi)\geq 0$. Thus, $\calR^*$ is indeed a dual dissipation
potential.

To derive the induced gradient-flow evolution for the GS
$(\PROB(S),\calE,\calR^*)$ we observe
\[
\rmD_\xi\calR^*(\rho,{-}\rmD\calE(\rho))[\eta] =
\rmD_\xi\calH(\rho,0)[\eta] = 
\int_S \Big(\ee^0({-}\eta)\bbQ (\ee^0)+\ee^0 \bbQ(\ee^0\eta) \Big)
\dd\rho =\int_S \eta \dd(\bbQ^*\rho) ,
\]
where we used $\bbQ1\equiv 0$. This provides
$\dot\rho=\rmD_\xi\calR^*(\rho,{-}\rmD\calE(\rho)) = \bbQ^*\rho$,
which is the expected Kolmogorov forward equation.

\subsubsection{A finite-state Markov process}
\label{sss:FiniteMP}

We consider the finite state space $S=\{1,\ldots,I\}$ such that 
\[
\PROB(S)=\bigset{ \rho=\bfc=(c_1,\ldots,c_I)\in [0,1]^I}{ \ts\sum_{i=1}^I c_i=1}.
\]
The Kolmogorov forward equation is the ODE 
\[
\dot\bfc = \bbA \bfc \quad \text{with }\bbA\in \R^{I\ti I}.
\]
Note that the Markov generator is given by $\bbQ_\mafo{finite}=\bbA^\top$, and the
conditions for a Markov generator are 
\[
\bbA_{ij}\geq 0 \text{ for all } i\neq j \quad \text{and}\quad 
\forall\,i=1,\ldots,I:\ 0=\sum_{j=1}^I \bbA_{ji}.
\]
We further assume that there is a unique positive steady state
$\pi=\bfw\in \PROB(S)$ such that the detailed-balance condition holds,
namely $\bbA_{ij}w_j=\bbA_{ji}w_i$. 

Thus, the induced energy functional is $\calE_{\mafo{finite}}(\bfc)=\frac12
\sum_{i=1}^I c_i\lambda_\rmB(c_i/w_i)$. To calculate the dissipation
potential we use that $\bbQ$ can be split 
\[
\bbQ_{\mafo{finite}}=\sum_{i=1}^{I-1}
\sum_{j=i+1}^I\bbQ^{i\leftrightarrow j} \text{ with }\bbQ^{i\leftrightarrow j}:= m_{ij}\Big(
\frac1{w_j} \bfe^j{\otimes}\bfe^i + \frac1{w_i}
\bfe^i{\otimes}\bfe^j\Big),
\]
where $m_{ij}:=\bbA_{ij}w_j=\bbA_{ji}w_i$ and $\bfe^k$ denotes the
$k$th unit vector in $\R^I$. Using the linearity in $\bbQ_{\mafo{finite}}$ of the
formula \eqref{eq:calR-LDP} for $\calR_{\mafo{finite}}$ we can first calculate
\begin{align*}
\calR^*_{i\leftrightarrow j}(\bfc,\bfxi)&= \sum_{k,l=1}^I\bigg[
\big(\frac{c_k}{w_k}\big)^{1/2}\ee^{-\xi_k} \bbQ^{i\leftrightarrow
  j}_{kl} \big(\frac{c_l}{w_l}\big)^{1/2}\ee^{\xi_l} - 
\big(\frac{c_k}{w_k}\big)^{1/2}
 \bbQ^{i\leftrightarrow j}_{kl} \big(\frac{c_l}{w_l}\big)^{1/2} \bigg]w_k
\\
&= m_{ij} \Big(\frac{c_ic_j}{w_iw_j}\Big)^{1/2} \big(
\ee^{\xi_i-\xi_j} + \ee^{\xi_j-\xi_i} -2 \big). 
 \end{align*}
Summing these terms and using the function $\CCC^*$ we find
\[
\calR_{\mafo{finite}}^*(\bfc,\bfxi)=\frac12\sum_{i=1}^{I-1} \sum_{j=i+1}^I
m_{ij} \Big(\frac{c_ic_j}{w_iw_j}\Big)^{1/2}\: 
\CCC^*\big(2(\xi_i{-}\xi_j)\big) 
\]
and conclude by Theorem \ref{th:LDgivesGS} that the equation
$\dot\bfc=\bbQ_{\mafo{finite}}^\top\bfc$ is induced by the GS
$(\PROB(S),\calE_{\mafo{finite}},\calR_{\mafo{finite}})$.

\subsubsection{Linear reaction-diffusion systems}
\label{sss:LinRDS}

We now return to RDS as discussed in
Section \ref{sss:nonlRDS}, but now consider only linear reactions where
all stoichiometric vectors $\bfalpha^r$ and $\bfbeta^r$ are given by
unit vectors $\bfe^i$ and $\bfe^j$, respectively. This means that the
reaction is a simple exchange reaction $C_i \rightleftharpoons C_j$. 
The linear RDS on a bounded smooth domain $\Omega \subset \R^d$ takes the form 
\begin{equation}
  \label{eq:lRDS}
  \dot \bfc = \bbD \Delta \bfc + \bbA \bfc, \quad \text{where } 
   \bbD=\diag(\delta_j)_{j=1,\ldots,I} \text{ with }\delta_j \geq 0
\end{equation}
complemented by no-flux boundary conditions. The matrix
$\bbA$ is as before. 
Now, $c_j(t,\cdot) \in \rmL^1(\Omega)$ is the nonnegative
concentration of the chemical species $C_i$. 

This system can be understood as the Kolmogorov forward equation on
the state space $S=\Omega\ti \{1,\ldots,I\}$, where the random variable
$Y_t=(X_t,i(t))$ undergoes a Brownian motion in $\Omega$ with
diffusion constant $\delta_j$ as long as $i(t)=j$. At discrete times
the particle can change its type within $\{1,\ldots,I\}$, according to the
jump process induced by the generator $\bbQ_{\mafo{finite}}=\bbA^\top$,
and then continue a Brownian motion with the new diffusion
constant. The full generator is 
\[
(\bbQ f)(x,i)= \delta_i \Delta f(x,i) + \sum_{k=1}^I \bbA_{ki} f(x,k),
\quad  \nabla f(x,i)\cdot \nu=0 \text{ on }\pl\Omega.
\]

We now assume that the linear reaction system satisfies the 
detailed-balance condition, i.e.\ we assume that there is an 
equilibrium state $\bfw$ with $w_i>0$  and $\bbA_{ij} w_j= \bbA_{ji}w_i$ for
all $i$ and $j$. Then, the steady state $\pi \in\PROB(S)$ is given by the
product of the $d$-dimensional Lebesgue measure on $\Omega$ and
$\bfw$, up to a suitable normalization factor:
\[
\pi = \frac1Z\: \rmd x \otimes w 
\text{ where }Z=\sum_{i=1}^I w_i \:\mafo{vol}(\Omega). 
\]   
By normalizing
$\bfw$ suitably, we may assume $Z=1$ subsequently.  

Using the Neumann boundary conditions, it is easy to check that the
generator $\bbQ$ satisfies the detailed-balance condition
\eqref{eq:DBC-Markov} with respect to $\pi$.

Hence, we can apply Theorem \ref{th:LDgivesGS} which provides the
large-deviation GS for \eqref{eq:lRDS}. Note that
$\rho\in \PROB(S)$ is absolutely continuous with respect to $\pi$ if
and only if $\rho=(c_i\dd x)_{i=1,\ldots,I}$ with $\bfc=(c_i)_{i=1,\ldots,I} \in
\rmL^1_{\geq 0}(\Omega)^I$. Moreover,  $\frac{\rmd\rho}{\rmd \pi} =
\big( c_i/w_i\big)_{i=1,\ldots,I}$ shows that the probability density
$\frac{\rmd\rho}{\rmd \pi}$ equals the vector of relative
concentrations $c_i/w_i$. 

The driving functional $\calE$ is the relative entropy 
up to a factor 1/2,
viz.\
\begin{equation}
  \label{eq:lRDS-F}
  \calE(\rho)= \frac12 \int_S \lambda_\rmB\big(\frac{\rmd\rho}{\rmd
    \pi}\big) \dd\pi = \frac12\int_\Omega  \Big(  \sum_{i=1}^I\lambda_\rmB(
  c_i(x)/w_i) w_i \Big)\dd x. 
\end{equation}  
For calculating the dissipation potential $\calR$ we can take
advantage of the linearity in $\bbQ$. In fact, $\bbQ$ can be split
into $I$ diffusion processes and the reaction part, namely
$\bbQ=\sum_{i=1}^I\bbQ^{(i)}_{\mafo{diff}}+ \bbQ_{\mafo{finite}}$ with
$\bbQ_{\mafo{finite}}$ as above.  The corresponding functionals
$\calH_{j}$ for diffusion processes $\bbQ^{(j)}_{\mafo{diff}}$ take
the form
\begin{align*}
\calH_{j}(\rho,\xi)&= \delta_j \int_\Omega\Big( |\nabla \xi(x,j)|^2 + \Delta
\xi(x,j) \Big)\dd \rho(x,j).
\end{align*}  

The dual dissipation potential $\calR^*$ is obtained by replacing
$\xi(x,j)$ by $\xi(x,j){+} \frac12\log\big(\frac{c_j(x)}{w_j}\big)$
and $\rmd\rho(x,j)=c_j(x)\dd x$, where we also use
\[
\int_\Omega \nabla c_j(x)\cdot \nabla \xi(x,j)+ c_j(x)\Delta \xi(x,j)
\dd x = 
\int_{\pl\Omega } c_j(x) \nabla \xi(x,j)\cdot \nu \dd a = 0.
\]
Subtracting the term at $\xi=0$, writing $\bfxi=(\xi_j)_j$ with
$\xi_j(x)=\xi(x,j)$, and using the result for $\bbQ_{\mafo{finite}}$
from above, we arrive at the formula 
\begin{align*}
\calR^*(\bfc , \bfxi)& = \int_\Omega \sum_{j=1}^I  \delta_j
|\nabla\xi_j(x)|^2  c_j(x) \dd x   \\
& \quad+ \int_\Omega \sum_{i=1}^{I-1}
    \sum_{k=2}^I \frac{1}2 \,\sqrt{\bbA_{ki}c_i(x)\bbA_{ik}c_k(x)}  \:
     \CCC^*\big(2(\xi_k(x){-}\xi_i(x))\big)\dd x.
\end{align*}
This is the same dual
dissipation potential as given in \eqref{eq:nRDS-R*CCC}, except for
the factors $\frac12$ and $2$ outside and inside of $\CCC^*$. 
However, these scaling factors arise since in the large-deviation
result in Theorem \ref{th:LDgivesGS} a factor
$\frac12 $ appears in the definition of $\calE(\rho)$ 
(see \eqref{eq:lRDS-F}), which is one-half of the usual relative entropy.

\subsubsection{Large deviations for a membrane model}
\label{sss:MP-Memb}

We consider a diffusion equation in the interval $\Omega={]{-}1,1[}$,
where at $x=0$ there is a membrane giving rise to a transmission
condition. The Kolmogorov forward equation takes the form (where
$\dot{\ }=\pl_t$ and $'=\pl_x$)
\begin{align*}
\dot\rho&= a_\pm \rho'' \text{ for } {\pm\, x} \in {]0,1[},\qquad  
  0= a_\pm \rho'(\pm 1),  \\
& a_+ \rho'(0^+) = b\big(\rho(0^+)- \rho( 0^-)\big)= 
 a_- \rho'(0^-). 
\end{align*}
The last relation means first that the mass flowing out of ${]{-}1,0[}$ has
to equal the flow into ${]0,1[}$, and second that this flow is
proportional to the difference of the densities. 

The invariant measure is $\pi=\frac12 \dd x$, and the Markov
generator $\bbQ_{\mafo{memb}}$ takes the form 
\begin{align*}
(\bbQ_{\mafo{memb}} f)(x)&= a_\pm f''(x) \text{ for }\pm x\in{]0,1[}, 
\qquad 0= a_\pm f'(\pm1),\\ 
& a_+ f'(0^+)=b\big(f(0^+){-}f(0^-)\big) =  a_-f'(0^-). 
\end{align*}
The functional $\calH_{\mafo{memb}}$ takes the form 
\begin{align*}
&\calH_{\mafo{memb}}(\rho,\zeta)=
\int_{{]{-}1,0[}} a_- \big(\zeta''+(\zeta')^2\big)\rho\dd x + 
\int_{{]0,1[}} a_+ \big(\zeta''+(\zeta')^2\big)\rho\dd x \\ 
& \text{ with } \zeta'(\pm1)=0 \text{ and }a_+\zeta'(0^+)\ee^{\zeta(0^+)}=
b\big(\ee^{\zeta(0^+)}-\ee^{\zeta(0^-)}\big) = a_-\zeta'(0^-)\ee^{\zeta(0^-)}. 
\end{align*}
Inserting $\zeta=\xi+\frac12\log(2\rho)$ and doing an integration by
parts using the nonlinear boundary conditions one obtains the dual
dissipation potential 
\[
\calR_{\mafo{memb}} (\rho,\xi)=\int_{{]{-}1,0[}}\!\! a_-(\xi')^2\rho\dd x
+ \sqrt{\rho(0^-)\rho(0^+)}\,\CCC^*\big(2(\xi(0^+){-}\xi(0^-))\big) 
+\int_{{]0,1[}} \!\! a_+ (\xi')^2\rho\dd x, 
\]
which again features the non-quadratic dissipation function $\CCC^*$. 


\section{Evolutionary $\Gamma$-convergence}
\label{ss:EGC}

Following the notions in the survey 
\cite{Miel14?EGCG} we consider families 
of GS $(\bfX,\calE_\eps,\calR_\eps)_{\eps\in {]0,1[}}$ and ask the question
whether the solutions $u_\eps$ for these systems have a limit $u$ for
$\eps\to 0$ and whether $u$ is again a solution to a
GS $(\bfX,\calE_0,\calR_0)$. Ideally, one
might hope that it is sufficient for $\calE_\eps$ and $\calR_\eps$
to converge in a suitable topology to $\calE_0$ and $\calR_0$,
respectively. Such results indeed exist and can be found in the surveys
\cite{Serf11GCGF,Miel14?EGCG}. However, the aim of this work is to
highlight the fact that starting with  classical (i.e.\ quadratic) 
dissipation potentials $\calR_\eps$ we may end 
up with a limiting dissipation $\calR_0$ that is non-quadratic. Thus,
limits of classical GS may be generalized GS. First such examples were
given in \cite{Miel12ERID,MieTru12DVEC} in the context of plasticity. 

\subsection{pE-convergence of gradient systems}
\label{ss:pE-Conv}
We first recall the general definition of pE-convergence, which
is a short name for evolutionary $\Gamma$-convergence with
well-prepared initial conditions.  Hence, the letter``E'' stands for
both, `E'volutionary convergence and `E'nergy convergence, while the
letter ``p'' stands for well`P'reparedness of the initial conditions, i.e., 
$ \calE_\eps(0,u_\eps(0)) \to \calE_0(0,u^0){<}\infty $. 

\begin{definition}[pE-convergence of $(\bfX,\calE_\eps,\calR_\eps)$]
\label{def:pE-conv}  
We say that the generalized gradient systems
$(\bfX,\calE_\eps,\calR_\eps)$ pE-converge  to   
$(\bfX,\calE_0,\calR_0)$, and write
$(\bfX,\calE_\eps,\calR_\eps) \pEweak (\bfX,\calE_0,\calR_0)$,  if 
\begin{align}
   \label{eq:Def-pE}
\left. \ba{@{\!}c@{\!}} u_\eps:[0,T]\to \bfX\quad\mbox{}\\
  \text{ is sol.\,of\/ }(\bfX,\calE_\eps,\calR_\eps),\\
    u_\eps(0)\weak u^0, \text{ and}\\
 \calE_\eps(0,u_\eps(0)) \to \calE_0(0,u^0){<}\infty \ea \right\} 
 \:  \Longrightarrow  \: 
 \left\{ \ba{@{\!}c@{\!}} \exists\, u \text{ sol.\ of\/ }
   (\bfX,\calE_0,\calR_0) \text{ with  } u(0){=}u^0\\ 
 \text{ and a subsequence } \eps_k\to 0: \\
 \forall\, t\in {]0,T]}{:} \ u_{\eps_k}(t)\weak u(t) \text{ and}\\
 \mbox{}\qquad \qquad \  \calE_{\eps_k}(u_{\eps_k}(t)) \to \calE_0(u(t)).
 \ea \right.
\end{align}
\end{definition} 
Here $u_\eps \weak u$ means the weak convergence in the Banach space
$\bfX$.  We emphasize that the notion of pE-convergence asks for
convergence of both, the solutions and the energies, but not of the
dissipation potentials. However, using the EDP and the convergence of
the energies, we easily obtain convergence of the integrated
dissipations, namely
\begin{align*}
&\int_0^T \calR_{\eps_k}(u_{\eps_k}, \dot u_{\eps_k}){+} \calR^*_{\eps_k} 
\big(u_{\eps_k} , {-}\rmD\calE_{\eps_k}(u_{\eps_k})\big)\dd t \ = \ 
\calE_{\eps_k}(u_{\eps_k}(0))-\calE_{\eps_k}(u_{\eps_k}(T))
\\
&\to \ \calE_{0}(u(0))-\calE_{0}(u(T))\ = \
\int_0^T \calR_0(u,\dot u){+}\calR^*_0\big(u,{-}\rmD\calE_0(u)\big) \dd t.
\end{align*}

A first systematic study of evolutionary $\Gamma$-convergence relying
on gradient structures was initiated in Sandier-Serfaty
\cite{SanSer04GCGF}, see also \cite{Serf11GCGF,Miel14?EGCG}. In this
approach one derives \emph{sufficient conditions} for pE-convergence
based on a limiting passage in the EDB
\begin{equation}
  \label{eq:EDBeps}
  \calE_\eps(u_\eps(T)) + \int_0^T \calR_\eps(u_\eps,\dot u_\eps) 
+ \calR^*(u_\eps,{-}\rmD\calE_\eps(u_\eps) ) \dd t = \calE_\eps(u_\eps(0)). 
\end{equation}

We observe that on the right-hand side we have the initial
energy, which converges to the desired limit because of the
well-prepared initial conditions.  Thus, to obtain an (EDE) for the
limiting process it suffices to show three liminf estimates for the
terms on the left-hand side, namely
\begin{subequations}
  \label{eq:SanSer}
\begin{align}
  \label{eq:SanSer1}
&\liminf_{\eps\to 0} \calE_\eps(u_\eps(T)) \geq \calE_0(u(T));\\
  \label{eq:SanSer2}
&\liminf_{\eps\to 0} \int_0^T \calR_\eps(u_\eps,\dot u_\eps) \dd t  \geq 
  \int_0^T \calR_0(u,\dot u)\dd t;\\
  \label{eq:SanSer3}
&\liminf_{\eps\to 0} \int_0^T\calR^*_\eps\big(u_\eps,{-}\rmD\calE_\eps(u_\eps)\big) 
 \dd t \geq  \int_0^T\calR^*_0\big(u,{-}\rmD\calE_0(u)\big)\dd t.  
\end{align}
\end{subequations}
Of course, it is sufficient that these convergences hold only along (a
subsequence of) the solutions $u_\eps$ of (EDP). In the following
subsection, we  will generalize this approach by keeping the terms
$\calR_\eps$ and $\calR_\eps^*$ together.

\subsection{EDP-convergence for gradient
  systems}\label{ss:ECgGS}

Here we define a new notion of evolutionary $\Gamma$-convergence for
GS that, on the one hand, is more restrictive but, one the other
hand, gives a more precise information on the limiting dissipation
potential $\calR_0$. We use the
fact that we do not need to have the two convergences
\eqref{eq:SanSer2} and \eqref{eq:SanSer3} separately. Indeed, it is
sufficient that only the integral over the sum of the two terms in
\eqref{eq:EDBeps} converges. 
This approach relaxes the sufficient conditions for pE-convergence 
substantially, since
in the limit $\eps\to 0$ the different 
parts of the dissipation may be distributed differently. 

In the quadratic case we always have equidistribution
$\calR_\eps(u_\eps,\dot u_\eps)= \calR^*_\eps(u_\eps,{-}
\rmD\calE_\eps(u_\eps))$ for solutions $u_\eps$. So, if \eqref{eq:SanSer2} and
\eqref{eq:SanSer3} hold, we will still have equidistribution in the
limit, but only along the limit solutions, while the limit functionals
need not be quadratic. In \cite{Miel12ERID} it is shown that the limit
of a classical GS can be a rate-independent system, where we always
have $\calR^*_0(u,{-}\rmD\calE_0(u))=0$, so there is not even
equidistribution along solutions.

We also add a strengthening condition to obtain our notion of ``EDP
convergence'' by asking for the convergence not only along solutions,
but rather along a suitable general class of functions $u:[0,T]\to
\bfX$. For this, we associate with the GS
$(\bfX,\calE_\eps,\calR_\eps)_{\eps\in [0,1]}$ De~Giorgi's
dissipation functionals ${\mathscr D}_\eps$, which are defined as
\[
{\mathscr D}_\eps(u) := \int_0^T \calR_\eps(u,\dot u)+ \calR_\eps^*\big(u,{-}\rmD
\calE_\eps(u)\big)  \dd t  .
\]

\begin{definition}[EDP convergence]\label{def:EDP-conv}
  We say that the family $(\bfX,\calE_\eps,\calR_\eps)_{\eps>0}$
  converges in the EDP sense to $(\bfX,\calE_0,\calR_0)$, and shortly
  write $(\bfX,\calE_\eps,\calR_\eps) \ \EDPweak \
  (\bfX,\calE_0,\calR_0)$ if
\begin{subequations}
  \label{eq:defEDP}
\begin{align}
  \label{eq:defEDP.a}
& (\bfX,\calE_\eps,\calR_\eps) \ \pEweak \
  (\bfX,\calE_0,\calR_0), \\[0.3em]
  \label{eq:defEDP.b}
&\calE_\eps \Gweak \calE_0 \text{ in } \bfX, \text{ and }\\[0.3em]
  \label{eq:defEDP.c}
&\ba{@{}l} \wt u_\eps(\cdot) \overset{*}{\rightharpoonup} \wt u(\cdot) 
 \text{ in } \rmL^\infty([0,T];\bfX) \text{ and}\\
\sup_{\eps\in {]0,1]},\:t\in [0,T]} \calE_\eps(\wt u_\eps(t)) \leq C < \infty
\ea\!\Big\} \  \Longrightarrow \
 \liminf_{\eps\to 0} {\mathscr D}_\eps(\wt u_\eps) \geq
 {\mathscr D}_0(\wt u) . 
\end{align}
\end{subequations}
\end{definition}
We emphasize that in condition \eqref{eq:defEDP.c}, the functions $\wt
u_\eps$ are arbitrary and need not be solutions of the GS
$(\bfX,\calE_\eps,\calR_\eps)$. 
From this definition we see that the convergence conditions 
\eqref{eq:SanSer} obviously imply EDP convergence. However, we will
study cases, where \eqref{eq:SanSer} does not hold, but we still have
EDP convergence. 

We emphasize that EDP-convergence is to be expected whenever one uses
the EDP principle for establishing pE-convergence. Indeed, from
general arguments one presumes that De Giorgi's dissipation functional
${\mathscr D}_\eps$ has (after extraction of a subsequence) a
$\Gamma$-limit ${\mathscr D}_0$ in the form
\[
{\mathscr D}_0(u) = \int_0^T\!\!\!\calM_0(u(t),\dot u(t)) \dd t.
\]
From the lower semicontinuity of $\Gamma$-limits, one expects that
$\calM(u,\cdot)$ is convex. Hence, one can define $\calR_\calM$  via
\[
\calR_\calM(u,v):= \calM_0(u,v)- \calM_0(u,0)
\]
and hope that it is a dissipation potential. For this, one needs to
show (i) the positivity $\calR_\calM(u,v)\geq 0$ and (ii) $\calM_0(u,0)\geq
\calR_\calM^*(u,{-}\rmD\calE_0(u))$. Often, the positivity (i) follows
simply from the evenness $\calM_0(u,{-}v)=\calM_0(u,v)$ and convexity.  
Moreover, if it is
possible to show $\calM_0(u,v) \geq -\langle \rmD \calE_0(u), v
\rangle$ (this holds for $\eps>0$ in the form
$\calR_\eps(u,v){+}\calR^*_\eps(u,{-}\rmD\calE_\eps(u)) \geq -\langle
\rmD\calE_\eps(u),v\rangle$), then we find (ii) via the estimate
\[
\calR_\calM^*(u,{-}\rmD\calE_0(u)) =\sup_{v\in \bfX}\Big(
\langle{-}\rmD\calE_0(u), v\rangle -\calM_0(u,v)+ \calM_0(u,0)\Big)
\leq \calM_0(u,0).   
\]
Thus, we arrive at the desired EDE $\calE_0(u(T)) + \int_0^T \calR_\calM+
\calR_\calM^* \dd t \leq \calE_0(u(0))$ and EDP-convergence to 
$(\bfX,\calE_0,\calR_\calM)$ is established.  

It would be interesting to study more generally the relations between
pE and EDP-convergence. Obviously, showing the liminf
estimate for ${\mathscr D}_\eps$ is the major step in establishing
pE-convergence. Hence, it seems redundant to ask for the
pE-convergence explicitly, yet it is not obvious under what additional
condition (e.g.\ the validity of a suitable chain rule) we really can
deduce the pE-convergence from the liminf estimate for ${\mathscr D}_\eps$.
 
We end this section with two examples concerning  EDP-convergence. 
Example \ref{ex:pEnotEDP} shows that the model discussed in
\cite{Miel12ERID} satisfies pE-convergence but not EDP-con\-ver\-gence. 
Example \ref{ex:Tartar} emphasizes the fact that 
pE and EDP-convergence are not properties of an evolution equation $\dot u =
V_\eps(u)$ but of a GS $(\bfX,\calE_\eps,\calR_\eps)$. Indeed, for a
given equation one may 
have different gradient structures leading to different limits in the
EDP sense, which in turn generate different limit evolutions.

\begin{example}[pE-convergence without
  EDP-convergence]\slshape\label{ex:pEnotEDP}
We consider the wig\-gly-energy model introduced in \cite{AbChJa96KMWE}.
It is given via
the time-dependent GS $(\R,\calE_\eps,\calR_\eps)$ with 
\[
\calE_\eps(t,u)=\frac12 u^2 -\ell(t) u +r\eps\sin(u/\eps) \quad
\text{and} \quad \calR_\eps(\dot u)=\frac{\eps}2 \,\dot u^2.
\]
For sufficiently smooth loading curves $\ell:[0,T]\to \R$ it was shown
in \cite[Thm.\,3.2]{Miel12ERID} that the GS $(\R,\calE_\eps,\calR_\eps)$
pE-converge to the generalized GS
$(\R,\calE_{\mafo{play}},\calR_{\mafo{play}})$ defined in
\eqref{eq:Play}. Obviously, we have the uniform convergence
$\calE_\eps\to \calE_{\mafo{play}}$, while 
${\mathscr D}_\eps$ converges to a limit ${\mathscr D}_0$
that cannot be written in terms of $\calR_\mafo{play}+\calR^*_\mafo{play}$,
see \cite[Prop.\,3.1]{Miel12ERID}

\end{example}

\begin{example}[Different limit equations]\label{ex:Tartar}\slshape
Here, we provide
an example of an evolution equation $\dot u = V_\eps(u)$ with two
different gradient structures. Both gradient structures have an evolutionary
$\Gamma$-limit in the EDP sense, and the surprising fact is that the
generated limit evolutions are different. Thus, EDP-convergence and
pE-convergence are not properties of the family of evolution equations
$\dot u = V_\eps(u)$, but of the chosen gradient structures.

Consider $\Omega=[0,1]$ and let $\bfX=\bfM_{\geq0}(\Omega)$, the set
of nonnegative finite Radon measures. Moreover, consider a continuous
periodic function $\bbA:\R\to {]0,\infty[}$ such that
\[
0 < a_\mafo{min}:=\min\nolimits_{y\in \R} \bbA(y) < 
 \max\nolimits_{y\in \R} \bbA(y) =:a_\mafo{max} < \infty.
\]
With $a_\eps(x)=\bbA(x/\eps)$ we define the simple PDE 
\[
\dot u(t,x) = - a_\eps(x) u(t,x), \quad u(0,x)=u_0(x)>0\quad \text{ for } 
x \in \Omega. 
\]
We introduce two different gradient structures $(\bfX,\calE_\eps,\calR_\eps)$  and 
 $(\bfX,\wt\calE_\eps,\wt\calR_\eps)$ via
\[
\calE_\eps(u)=\int_\Omega a_\eps\dd u, \
\calR^*_\eps (u,\xi) = \int_\Omega  \frac12\xi^2 \dd u, \quad 
\wh\calE_\eps(u)=\int_\Omega \frac1{a_\eps}\dd u, \
\calR^*_\eps (u,\xi) = \int_\Omega  \frac{a_\eps^2}2\xi^2 \dd u.
\]
It is shown in \cite[Cor.\,3.8]{Miel14?EGCG} that these 
GS converge in the EDP sense to the limit systems $(\bfX,\calE_0,\calR_0)$  and 
 $(\bfX,\wt\calE_0,\wt\calR_0)$, respectively, where
\[
\calE_0(u)=\int_\Omega a_\mafo{min}\dd u, \
\calR^*_0 (u,\xi) = \int_\Omega  \frac12\xi^2 \dd u, \quad 
\wh\calE_0(u)=\int_\Omega \frac1{a_\mafo{max}}\dd u, \
\calR^*_0 (u,\xi) = \int_\Omega  \frac{a_\mafo{max}^2}2\xi^2 \dd u.
\]
In particular, the limit evolution for the first is $\dot u= -
a_\mafo{min} u$, while it is $\dot u=-a_\mafo{max} u$ for the second.
This is not a contradiction, but has its origin in the
well-preparedness condition for the initial data. No sequence can be
well-prepared for both systems, i.e.\ if $u_\eps \weaks u$ and 
$\calE_\eps(u_\eps)\to \calE_0(u)$, then we have
$\wh \calE_0(u)\lneqq \liminf_{\eps\to 0}\wh\calE_\eps(u_\eps)$, and
vice versa. 
\end{example}

\subsection{EDP-convergence for an ODE example}
\label{ss:ODE} 

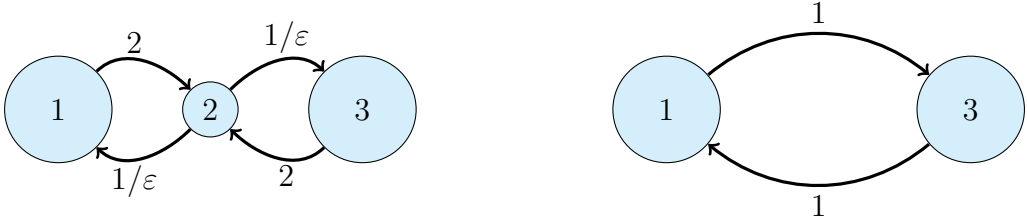
\begin{figure}
\centering
\begin{tikzpicture}
\tikzstyle{every node}=[draw,shape=circle,fill=cyan!15];
\node (v1) at (2,0) {$\quad1\quad$};
\node (v2) at (4,0) {$2$};
\node (v3) at (6,0) {$\quad3\quad$};
\node (vv1) at (10,0) {$\quad1\quad$};
\node (vv3) at (14,0) {$\quad3\quad$};
\tikzstyle{every node}=[];

\draw[->, very thick] (v1) to [out=45, in=135] (v2); 
\draw[->, very thick] (v2) to [out=45, in=135] (v3); 
\draw[->, very thick] (v2) to [out=-135, in=-45] (v1); 
\draw[->, very thick] (v3) to [out=-135, in=-45] (v2); 

\draw[->, very thick] (vv1) to [out=40, in=140] (vv3); 
\draw[->, very thick] (vv3) to [out=220, in=-40] (vv1);

\node[above] at (3,.6) {$2$}; 
\node[above] at (5,.6) {$1/\eps$}; 
\node[below] at (3,-.6) {$1/\eps$}; 
\node[below] at (5,-.6) {$2$}; 

\node[above] at (12,1) {$1$}; 
\node[below] at (12,-1) {$1$}; 
\end{tikzpicture}
\caption{Left: Three-state Markov process with high rate of leaving state
2. Right: The limit for $\eps\to 0$ gives a two-state Markov process.}
\label{fig:Markv3-2} 
\end{figure}

We discuss a very simple example of a discrete Markov process with
state space $S=\{1,2,3\}$. The jump rates are such that in the limit
$\eps\to 0$ the particles never stay in the state $2$. Thus, the
limiting Markov process has the state space $\{1,3\}$ only, see Figure
\ref{fig:Markv3-2}. We will
start with three different GS, namely (i) the quadratic one, where
both  $\calE_\eps$ and $\calR_\eps$ are quadratic, (ii) the entropic
one with classical $\calR^*_\eps$, and (iii) the entropic one with the
dual dissipation potential defined in terms of $\CCC^*$. 
The interesting point is that in the cases (i) and (iii) the limiting
GS obtained via EDP-convergence will still be in the same modeling
class. However, in case (ii) we will lose the classical GS and obtain
a generalized GS that cannot be described via $\CCC^*$. 

We consider the  Kolmogorov forward equation (here an ODE) of a
 Markov process on the state space $S=\{1,2,3\}$ given by 
\begin{equation}
  \label{eq:ODE01}
  \dot u = (2{+}\eps) \bma{ccc} -1&1/\eps&0\\ 1&-2/\eps&1\\
  0&1/\eps&-1\ema u,\quad \bfX:=\PROB(\{1,2,3\}). 
\end{equation}
The unique equilibrium  $w^\eps =\frac1{2{+}\eps} (1,\eps,1)^\top$
satisfies the detailed-balance condition \eqref{eq:DBC-Markov}.

The limit dynamics is easily obtained by setting $u_2=\eps r$ giving 
\[
\dot u_1= (2{+}\eps)(r{-}u_1),\quad 
\eps \dot r= (2{+}\eps)(u_1 {-}2r{+}u_3),\quad 
\dot u_3= (2{+}\eps)(r{-}u_3).
\]
Thus, in the limit $\eps \to 0$ we find $r=(u_1{+}u_3)/2$ and 
\[
\dot u_1= u_3-u_1, \quad 0= u_1-2r+u_3,\quad \dot u_3 = u_1-u_3.
\]
More precisely, if the initial condition satisfies 
$u^\eps(0)\to (p_0,0,1{-}p_0)$, then for all $t>0$ we have 
\[
u^\eps(t)\to \big(p(t),0,1{-}p(t)\big)^\top \quad \text{where } 
 \dot p(t)=1-2p(t),  \   p(0)=p_0.
\]
We will study the limit $\eps\to 0$ in several gradient
structures. For general strictly convex and superlinear functions
$\phi$ and  $\psi$ we consider the GS 
$(\bfX,\calE_\eps,\calR_\eps)$ with $ \calE_\eps(u)=\sum_{i=1}^3
w^\eps_i \phi(u_i/w_i^\eps)$ (recall
$w^\eps=\frac1{2+\eps}(1,\eps,1)^\top$) and 
\begin{align*}
  &\calR^*_\eps(u,\xi) = \sum_{j=1}^2 a_j^\eps(u) \,
       \psi^*\big( \xi_{j+1}{-}\xi_j\big),
 \quad 
\text{where } a_j^\eps(u):=\frac{u_{j+1}/w_{j+1}^\eps -
      u_j/w_j^\eps}{(\psi^*)'\big(\phi'(\tfrac{u_{j+1}}{w_{j+1}^\eps}) 
    {-} \phi'(\tfrac{u_j}{w_j^\eps})\big)}.  
\end{align*} 
Using the fact that $v =\dot u$ satisfies $v_2= - v_1-v_3$ we
obtain the primal dissipation potential $\calR_\eps$ in the form 
\begin{equation}
  \label{eq:ODE.5}
  \calR_\eps(u,v)= a^\eps_1(u) \psi\big( \frac{v_1}{a^\eps_1(u)}\big) 
      + a^\eps_2(u) \psi\big( \frac{v_3}{a^\eps_2(u)}\big).
\end{equation}
Indeed, in general $\calR_\eps$ is an inf-convolution, since
$\calR_\eps^*$ is a sum over two terms. However, here we can eliminate
$v_2$ and argue as follows:
\begin{align*}
\calR_\eps(u,v)&=\sup\set{\xi\cdot v - \calR_\eps^*(u,\xi) }{
    \xi_1+\xi_2+\xi_3=0}  
\\ &=\sup\set{(\xi_1{-}\xi_2)v_1 + (\xi_3{-}\xi_2)v_3 -
  \calR_\eps^*(u,\xi)}{\xi_1+\xi_2+\xi_3=0 } ,
\end{align*}
which gives the desired result, since $\calR_\eps^*$ only depends on 
$\xi_1{-}\xi_2$ and $\xi_2{-}\xi_3$. 

We now state the result on EDP-convergence of $(\bfX,\calE_\eps,
\calR_\eps)$ to $(\bfX,\calE_0,\calR_0)$. Since the limiting GS can be
described on the much smaller set $\PROB(\{1,3\})$, which we
identify with $\bfY:=[0,1]$, we can formulate the limit GS in terms of
the reduced GS $([0,1],\bfE,\bfR)$.

\begin{theorem}\label{th:ODE-conv}
We have $(\bfX,\calE_\eps, \calR_\eps) \EDPweak (\bfX,\calE_0,\calR_0)$ with 
\begin{align*}
\calE_0(u)&=\left\{\ba{cl} \bfE(p)&
  \text{for }u=(p,0,1{-}p)^\top,\\
  \infty& \text{otherwise}; \ea \right. \text{ with } \bfE(p)=
\frac12\phi(2p)+\frac12\phi(2{-}2p)\\
\calR_0(u,\dot u)&=\left\{\ba{cl} \bfR(p,\dot p)&\text{for }(u,\dot
  u)=(p,0,1{-}p,\dot p,0,{-}\dot p)^\top, \\ \infty &\text{otherwise},
  \ea\right.
\end{align*}
where $\bfR$ is given via its Legendre dual 
\begin{align*}
\bfR^*(p,\eta)&:=\sigma(p) + \sup_{z>0}\left\{ \inf_{\tau\in \R}\Big(
    \wh a(p,z)\psi^*(\eta{-}\tau) {+} \wh a(1{-}p,z) \psi^*(\tau) \Big)
          -\Sigma(p,z) \right\},\\[0.5em]
\begin{split}
\text{where}\quad \Sigma(p,z)&:=\wh a(p,z)\psi^*\big(\phi'(2p){-} \phi'(z)\big) 
+ \wh a(1{-}p,z)\psi^*\big(\phi'(2{-}2p) {-} \phi'(z)\big), \\[0.5em]
\wh a(p,r)&:=(2p{-}r)/\big\{(\psi^*)'\big(
\phi'(2p){-}\phi'(r)\big)\big\},\text{ and} 
\\[0.5em]
 \sigma(p)&:= \inf\set{\Sigma(p,z)}{z>0}.
\end{split}
\end{align*}
\end{theorem}

In particular, this implies that 
the limiting ODE $\dot p= 1-2p$ is induced by the reduced
generalized GS $([0,1],\bfE,\bfR)$, i.e.\ $\dot p = 1-2p \ = \
\rmD_\eta \bfR(p, -\rmD\bfE(p))$.  
The above theorem follows directly from the next proposition and
the general theory described in Section \ref{ss:ECgGS}.  For the
energy-dissipation principle we consider De 
Giorgi's dissipation functional
\[
{\mathscr D}_\eps(u)=\int_0^T M_\eps (u(t),\dot u(t)) \dd t \quad \text{with } 
M_\eps(u,v)=\calR_\eps(u,v) + \calR_\eps^*(u,-\rmD\calE_\eps(u)).
\]

\begin{proposition}\label{pr:ODE-M0}
We have the $\Gamma $-limits
$\calE_\eps\overset\Gamma\to \calE_0$, $M_\eps \overset{\Gamma}\to M_0$, 
and  ${\mathscr D}_\eps \Gweak {\mathscr D}_0$  with 
\begin{align*}
&M_0(u,v)= \left\{\ba{cl} \bfR(p,\nu) {+}
  \bfR^*\big(p,\phi'(2p){-}\phi'(2{-}2p)\big) & \text{for 
  }(u,v)=(p,0,1{-}p,\nu,0,-\nu)^\top, \\ \infty
  &\text{otherwise};    \ea \right. 
\end{align*}
and ${\mathscr D}_0(u)=\int_0^T M_0(u(t),\dot u(t)) \dd t$, 
where $\bfR$ is given as in Theorem \ref{th:ODE-conv}. 
\end{proposition}
\begin{proof} The convergence $\calE_\eps \Gto \calE_0$ follows easily
  from the explicit form $\calE_\eps(u)=\sum_1^2 w_i^\eps
  \phi(u_i/w^\eps_i)$, the convergence $w^\eps_2\to 0$, and the
superlinearity of $\phi$. 

To simplify the $\Gamma$-limit of ${\mathscr D}_\eps$ we introduce the
scaling $u_2= \eps r$ and use that $\dot u_2 = \eps \dot r$
does not explicitly appear in $M_\eps$, see 
$\calR_\eps$ in \eqref{eq:ODE.5}.  However, the relation
$u_1+u_2+u_3\equiv 1$ now takes the $\eps$-dependent form $u_1+ \eps r
+ u_3\equiv 1$.  Moreover, defining $\wt M_\eps( u_1,r,u_3,\dot
u_1,\dot u_3)= M_\eps (u_1,\eps r,u_3,\dot u)$ shows that $\wt M_\eps$
is continuous in $\eps\in [0,1]$ and all the arguments.

Hence, when considering a sequence of functions $u^\eps:[0,T]\to \bfX$ with
$u^\eps \weaks u_0 $ and $ \calE_\eps(u^\eps(t)) \leq C< \infty$  
we find $u_0(t)=(p(t),0,1{-}p(t))^\top$ and 
\begin{align*}
&\liminf_{\eps\to 0} {\mathscr D}_\eps(u^\eps) \geq \int_0^T m(p(t),\dot p(t)) \dd t 
   \text{ with }\\
&m(p,v):=
\inf_{z>0}\Big\{  
\wh a(p,z)\psi\big(\frac{v}{\wh a(p,z)}\big) 
+ \wh a(1{-}p,z)\psi\big(\frac{v}{\wh a(1{-}p,z)}\big) + \Sigma(p,z) \Big\}.
\end{align*}
Note that we have no control over the variable $r(t)=u^\eps_2(t)/\eps$ in the
limit $\eps\to 0$. So, we simply minimize over all possible values $z\in
{[0,\infty[}$, which certainly provides a good lower bound. Moreover,
recovery sequences for the convergence ${\mathscr D}_\eps \Gweak {\mathscr D}_0$ can
be obtained in the form $u^\eps(t)=(p(t),0,1{-}p(t)) + \eps \zeta(t)
({-}p(t),1,p(t){-}1)$, where $\zeta(t)$ is the minimizer in the
definition of $m(p(t),\dot p(t))$. 

Obviously, the definition of $\sigma$ gives the relation $\sigma(p)=m(p,0)$. 
Thus, we have derived the reduced dissipation potential $\bfR$ in the 
form $\bfR(p,v)= m(p,v)-\sigma(p)$. Doing a Legendre transform with respect
to $v$ we obtain the form of $\bfR^*$ given in Theorem \ref{th:ODE-conv},
since the sum turns into an inf-convolution.
\end{proof}

Next we consider three different choices for $\phi$ and $\psi$.

\subsubsection{Quadratic energy and dissipation}
\label{sss:Quadratic} 

First, we consider the case   
\[
\phi(r)=\frac12\,r^2 \quad \text{and} \quad  \psi(v)=\frac12\,v^2,
\]
which gives $\bfE(p)=p^2+ (1{-}p)^2$ and $\wh a(p,r)\equiv 1$ and
simplifies all expressions considerably:
\[
\Sigma(p,z)=\frac12(2p{-}z)^2 + \frac12(2{-}2p{-}z)^2, \quad 
\sigma(p)=(1{-}2p)^2, \quad \bfR(p,\eta)=\frac14\,\eta^2. 
\]

Here, the crucial point in the definition of $\bfR^*$ is that the
inf-convolution involving $\psi^*$ and $\tau$ does not involve any
dependence on $z$, so that the term $\sigma(p)$ exactly cancels
$\sup_{z>0}{-}\Sigma(p,z)$, which is generally not the case.  Thus,
the limiting GS is $([0,1],\bfE,\bfR)$ where $\bfR(p,\nu)=\nu^2$ is
quadratic, and $([0,1],\bfE,\bfR)$ is again a classical GS.

\subsubsection{Entropic energy and $\CCC$-type dissipation}
\label{sss:EE+exp}

Next, we consider the case of the Boltzmann entropy and the 
dissipation defined in terms of $\psi=\CCC$, which coincides with 
Section \ref{ss:Markov} except for the
trivial scaling factor 2: 
\[
  \phi(r)=\lambda_\rmB(r)=r\log r - r +1 \quad \text{and} \quad  
  \psi^*(\xi)= \CCC^*(\xi)=4\big(\cosh(\xi/2) -1\big).
\]
This gives the reduced energy functional  $\bfE(p)=\frac12\lambda_\rmB(2p)+
\frac12 \lambda_\rmB(2{-}2p)$ and 
$\wh a(p,z)=  
\sqrt{2pz}$. In the latter expression and in the definition of
$\Sigma$ we profit from the interaction of 
$\lambda'_\rmB(r)=\log r$ and the exponential form of $\CCC^*$, viz.\ 
\[
\Sigma(p,z) = 2(\sqrt{2p}- \sqrt z)^2 + 2 (\sqrt{2{-}2p} - \sqrt z)^2 = 4 - 4
b_p \sqrt z + 4z\text{ with }b_p = \sqrt{2p} + \sqrt{2{-}2p}.
\]
Minimizing in $z>0$ we arrive at 
\[
\sigma(p)=4-b_p^2 = 2 -4 \sqrt{p(1{-}p)} = 2(\sqrt p - \sqrt{1{-}p})^2. 
\]

For calculating $\bfR^*$ we first observe, for $a,b>0$,  the formula 
\[
\inf_{\tau\in \R} \Big( a \CCC^*(\tau) + b \CCC^*(\xi{-}\tau)\Big)
\ = \ 4 \sqrt{(a{+}b)^2 + \frac{ab}2 \CCC^*(\xi)} - 4(a{+}b),  
\] 
which follows by writing the left-hand side via $\CCC^*(\xi{-}\tau)=
2\ee^\xi/x + 2 \ee^{-\xi}x -4$, where $x=\ee^\tau$, and minimizing in $x>0$.
With $a=\wh a(p,z)$ and $b=\wh a(1{-}p,z)$ we find 
\begin{align*}
\bfR^*(p,\eta)&= \sigma(p) + \sup_{z>0} \Big(
  4\sqrt z \sqrt{b_p^2 + \sqrt{p(1{-}p)} \CCC^*(\eta)} - 4b_p \sqrt z
  -\Sigma(p,z) \Big)\\
&= \sigma(p) + \sup_{z>0} \Big(
  4\sqrt z \sqrt{b_p^2 + \sqrt{p(1{-}p)} \CCC^*(\eta)}  - 4 - 4z \Big) \\
& = \sigma(p) -4+  b_p^2 + \sqrt{p(1{-}p)}\, \CCC^*(\eta)\ = \ 
  \sqrt{p(1{-}p)}\, \CCC^*(\eta). 
\end{align*}
We emphasize that in this minimization with respect to $z$ it is crucial to
keep the terms involving the dual dissipation potential $\CCC^*(\eta)$  
and the term $\Sigma$ together. 
 
We observe that the resulting gradient structure is again the 
structure, which is  obtained from the large-deviation
principle of Section \ref{ss:Markov}. 
This confirms the statement that gradient structures obtained
from the large-deviation theory are very stable against taking further
limits in the sense of EDP-convergence, see Figure \ref{fig:LDPvsEDP}

\subsubsection{Entropic energy and quadratic dissipation}
\label{sss:EE+quad}

In \cite{Maas11GFEF,ErbMaa12RCFM,Miel13GCRE,MaaMie15?GSRC1,MaaMie15?GSRC2}
the relative entropy was used for the energy functional
$\calE$ and a quadratic dissipation leading to a classical gradient
system:
%
%
\[  
\phi(r)=\lambda_\rmB(r)=r\log r - r +1 \quad \text{and} \quad 
\psi^*(\xi)=\frac12 \xi^2.
\]
We obtain the same limit energy $\bfE(p)=\frac12\lambda_\rmB(2p)+ \frac12
\lambda_\rmB(2{-}2p)$ as in the previous case, but the functions
$a_j^\eps(u)$ are quite different as they involve the logarithmic mean
$\Lambda(r,s)=\frac{r\ - \ s}{\log r{-}\log s}$. Indeed we have
$\wh a(p,z)= \Lambda(2p,z) = \frac{2p\  -\  z}{\log (2p){-}\log z}$.
We further obtain the functions 
\[
\Sigma(p,z)=\frac12\Big( (2p{-}z)\big(\log(2p)-\log z\big) + (2{-}2p{-}z)
\big( \log(2{-}2p)-\log z\big) \Big) 
\]
and have no explicit formula for $\sigma(p)=\inf_{z>0} \Sigma(p,z)$. 
In the definition of $\bfR^*$ we can do the inf-convolution explicitly, since
$\psi^*$ is quadratic, so we find the formula
\[
\bfR^*(p,\eta)= \sigma(p) + \sup_{z>0} \Big( \frac{\wh a(p,z)\wh a(1{-}p,z)}{
2(\wh a(p,z){+}\wh a(1{-}p,z))} \,\eta^2\;- \Sigma(p,z)  \Big).
\]
We claim that the growth of $\bfR^*(p,\eta)$ is no longer quadratic, but 
exponential. For this we insert $z=\ee^{b\eta}$ for $\eta\gg1$ 
for some $b\in\left]0,1/2\right[$
into the supremum to obtain a
lower bound. From  $\Sigma(p,z)\approx z \log z$ and $\wh a(p,z)
\approx z/\log z$ for $z\to \infty$ we find the asymptotic lower bound 
\[
\bfR^*(p, \eta) \ 
\raisebox{0.3em}{$>$}\hspace{-0.7em}\raisebox{-0.4em}{$\approx$} 
 \ \Big(\frac1{4b}-b\Big)\eta \ee^{b \eta}. 
\]
Hence, we see that the growth is at least as $\ee^{b|\eta|}$ for all $b \in
{]0,1/2[}$. Moreover, we expect that the function $\bfR^*(p,\eta)$ does not
have a product structure $b(p)\Psi(\eta)$ any more.

Thus, we see that the classical gradient structure for the relative entropy
is not stable under EDP, in general.  
Nevertheless, in
\cite{GigMaa13GHCD,DisLie13?GSMC,MaaMie15?GSRC1}  evolutionary
$\Gamma$-limits between discrete Markov processes and continuous
Fokker-Planck equation are studied, where the classical gradient
structure survives.

\section{The membrane as a thin-layer limit}
\label{se:Membrane}

In our first major application of the EDP-convergence as a microscopic
origin of generalized GS, we follow \cite{Lier12VME,Lier13PBBS} and 
consider a one-dimensional
diffusion equation with a thin layer of very small diffusivity. Assuming
that the diffusion coefficient and the width of the layer scale in the
proper way, we will arrive at a membrane model in the limit. While the
limit passage of the linear diffusion problem to the linear
transmission problem at the membrane can be done directly (or with the
quadratic gradient structure, see \cite{Lier13PBBS}), we prefer to do
the somewhat more elaborate EDP-limit using the GS with the relative
entropy as energy functional and the classical dissipation potential of
Wasserstein type. This case was already studied in
\cite[Sec.\,3.2]{Lier12VME} in a more special setting and without
explicitly calculating $\calR_0^*$.

We start from the equation 
\begin{equation}
  \label{eq:Memb-1}
  \dot u = \Big( a_\eps(x) \big( u' + u V_\eps'(x)\big)' \Big)
  \ \text{ in }\Omega:= {]{-}1,1[},\quad \pl_x u(t,\pm1) + u(t,\pm1)
  V'_\eps(\pm1)=0, 
\end{equation}
where $\dot{\ }=\pl_t$ and $'=\pl_x$.  By our choice of the boundary
conditions, the total mass $\int_\Omega u(t,x) \dd x =1$ is conserved,
thus we can interpret the equation as the Fokker-Planck equation of a
Markov process. Defining the equilibrium density
\[
w_\eps(x) = \frac1{Z_\eps} \ee^{-V_\eps (x)} \quad \text{with }Z_\eps=\int_{-1}^1
\ee^{-V_\eps(x)} \dd x,
\]
we have the GS $( \PROB(\Omega),\calE,\calR^*_\eps)$
with 
\[
\calE_\eps(u)=\int_\Omega \lambda_\rmB \big(u(x)/w_\eps(x)\big) w_\eps(x)\dd x \quad
\text{and} \quad \calR_\eps^*(u,\xi)= \frac12\int_\Omega a_\eps(x)
u(x) \xi'(x)^2 \dd x.
\]

The nontrivial behavior happens in the thin layer given by the
small interval $[0,\eps]$. In particular, we allow $a_\eps$ and $V_\eps$ to
depend non-trivially on $x$: We assume that there are functions
$a_*,a_+,V_*,V_+\in \rmC^1([0,1])$ and $a_-,V_-\in \rmC^1([-1,0])$ such that
$a_*(x),a_+(x),a_-(-x) \geq \underline a>0 $ for all $x\in [0,1]$,
$V_-(0)=V_*(0)$, $V_*(1)=V_+(0)$, and 
\begin{equation}
  \label{eq:Memb-3}
  a_\eps(x)=\left\{\ba{cl}a_+(x) &\text{for }x>\eps,\\
\eps a_*(x/\eps)& \text{for }x\in [0,\eps],\\   
  a_-(x)&\text{for }x<0,\ea \right. \quad \text{and} \quad 
  V_\eps(x)=\left\{\ba{cl}V_+(x{+}\eps) &\text{for }x>\eps,\\
V_*(x/\eps)& \text{for }x\in [0,\eps],\\   
  V_-(x)&\text{for }x<0.\ea \right. 
\end{equation}
Here $V_\eps$ is constructed to be continuous on
$\ol\Omega=[-1,1]$, while $a_\eps$ has jumps for $x\in \{0,\eps\}$.

The pE-convergence result established in
\cite[Sec.\,3]{Lier12VME} states that the limiting
system is given as a membrane problem, where the thin layer is
replaced by  a transmission condition. The interesting point is that
the EDP-convergence reveals that the limiting GS is no longer
classical but involves $\CCC$ for the jump of the driving
forces at the membrane. 

For passing to the limit we note that the function $w_\eps$ converges
pointwise to the limit
\begin{equation}
  \label{eq:Memb.w0}
  w_0(x)= \bigg\{ \ba{cl}\frac1{Z_0}\ee^{-V_+(x)}&\text{for }x>0,\\
\frac1{Z_0}\ee^{-V_-(x)}&\text{for }x<0,  \ea  \quad \text{with } 
Z_0= \int_{-1}^0\!\! \ee^{-V_-(x)}\dd x+ \int_0^1\!\! \ee^{-V_+(x)} \dd x,
\end{equation}
which may be discontinuous at $x=0$, but has well-defined limits 
$w_0(0^-)$ and $w_0(0^+)$ from the left and from the right,
respectively. This limit is totally independent of the potential $V_*$
inside the layer. The influence on the layer potential $V_*$ and the
layer diffusion profile $a_*$ will only survive in one coefficient
$A_*$.

\begin{theorem}[Membrane limit]\label{th:Membrane}
$(\PROB(\Omega),\calE_\eps,\calR^*_\eps) \ \EDPweak \
(\PROB(\Omega),\calE_0,\calR^*_0)$, where
\begin{align*}
\calE_0(u)&=\int_\Omega \lambda_\rmB (u/w_0) w_0\dd x \quad \text{and} \\
\calR^*_0(u,\xi)&= \int_{{]-1,0[}} \frac{a_-}2 (\xi')^2 u\dd x
 + \int_{{]0,1[}} \frac{a_+}2 (\xi')^2u \dd x \\ 
&\quad  + A_* \sqrt{\tfrac{u(0^-)u(0^+)}{w_0(0^-)w_0(0^+)}} 
\CCC^*\big( \xi(0^+){-}\xi(0^-)\big) \quad\text{where }
A_*= \Big( \int_0^1 \tfrac{Z_0 \ee^{V_*(y)}}{a_*(y)} \dd y \Big)^{-1}. 
\end{align*}
\end{theorem}

Before we go into the details of the proof, some comments are in
order.  First, we emphasize that the constant $Z_0$ in the definition
of the coupling coefficient $A_*$ is not related to $V_*$, 
but only depends on $V_\pm$, see \eqref{eq:Memb.w0}. Hence, for 
a large barrier $V_*$ the transmission coefficient $A_*$ becomes indeed small.

Second, the limiting equation is a PDE in the subdomains
$\Omega_-={]{-}1,0[}$ and $\Omega_+={]0,1[}$ coupled by a transmission
condition. It can be obtained easily by considering
test functions $\wh\xi \in \rmH^1(\Omega_-) \ti \rmH^1(\Omega_+)$ in the
weak form $ \int_\Omega \dot u \wh \xi \dd x = 
\rmD_\xi \calR^*_0\big(u,{-}\rmD \calE_0(u) \big)[ \wh\xi]$.
Using the fact that $\wh\xi$ may have a jump at $x=0$, the transmission
conditions arise via the boundary terms when integrating by
parts. We arrive at
\begin{align*}
\dot u&= \Big( a w_0 \big(u/w_0)\big)' \Big)' \quad \text{for }x 
\in \Omega_- \cup \Omega_+. \\
0&= a_+(0)w_0(0^+)\big(u/w_0)\big)'(0^+) - A_*\Big(\tfrac{u(0^+)}{w_0(0^+)}
                    - \tfrac{u(0^-)}{w_0(0^-)}\Big), \\
0&= a_- (0)w_0(0^-)\big(u/w_0)\big)'(0^-) -A_*\Big( \tfrac{u(0^+)}{w_0(0^+)}
                    - \tfrac{u(0^-)}{w_0(0^-)}\Big), \\
0&= a_\pm(x)w_0(x)\big(u/w_0)\big)'(x) \quad \text{at }x\in \{-1,1\}. 
\end{align*}
We refer to \cite{GliMie13GSSC} for a similar derivation of more general
nonlinear transmission conditions and active interface conditions
using gradient structures.  
 
Finally, we remark that the primal dissipation potential can be
written using the integration operator $I[\dot u](x):= \int_{-1}^x
\dot u(y)\dd y = -\int_{x}^1 \dot u(y) \dd y$, where
the last relation follows from $\int_{-1}^1 \!\dot u \dd y=0$,
which in turn is due to $u(t)\in \PROB(\Omega)$. Noting that the functions
$\xi$ may have a jump at $x=0$, one has the identity
\[
\int_{-1}^1 \xi \dot u\dd x = - \int_{-1}^0 I[\dot u] \xi'\dd x -
 I[\dot u](0)\big(\xi(0^+){-}\xi(0^-)\big) - \int_0^1 I[\dot u]\xi' \dd x.
\]
Here $I[\dot u](0)$ is the flux through the membrane, which is
thermodynamically conjugate to the jump $\xi(0^+){-}\xi(0^-)$ in the
driving forces.  With this and $I[\dot u](-1)=0=I[\dot u](1)$ the
evaluation of the Legendre transform for $\calR_0^*$ yields the primal
dissipation potential
\begin{equation}
  \label{eq:Memb.calR0}
\calR_0(u,\dot u) =  \int_{-1}^0  \frac{I[\dot u]^2}{2a_- u } \dd x 
 + \int_0^1 \frac{I[\dot u]^2}{2a_+ u} \dd x+
A_* \sqrt{\tfrac{u(0^-)u(0^+)}{w_0(0^-)w_0(0^+)}}
\CCC\Big( \frac{I[\dot u](0)}{A_* 
      \sqrt{\tfrac{u(0^-)u(0^+)}{w_0(0^-)w_0(0^+)}}} \Big).
\vspace{0.5em}
\end{equation}

\noindent{\itshape Proof of Theorem \ref{th:Membrane}:}
We first observe that $\calE_\eps \Gweaks \calE$ using
\cite[Lem.\,9.4.2]{AmGiSa05GFMS}. Moreover, pE-convergence was
established in \cite[Sec.\,3.2]{Lier12VME}. It remains to establish
the liminf estimate \eqref{eq:defEDP.c}, where De Giorgi's dissipation
functional ${\mathscr D}_\eps$ takes the explicit form
\begin{align*}
&{\mathscr D}_\eps(u):=\int_0^T\int_\Omega \frac1{2a_\eps u}I[\dot u]^2 +
\frac{a_\eps u}2 \Big(\big(\log(\tfrac{u}{w_\eps})\big)'\Big)^2 \dd x \dd t.
\end{align*}

\textit{Step 1. Blow up = transformation from ${\mathscr D}_\eps$ to
  $\wh{\mathscr D}_\eps$:}\/ 
To study the $\Gamma$-limit ${\mathscr D}_0$ of ${\mathscr D}_\eps$ we blow up the
thin layer such that its transformed thickness becomes of order
one. For this we use $Y_\eps:[-1,1]\to [-1,2]$ and its inverse $X_\eps
= Y_\eps^{-1} $:
\[
Y_\eps(x)=\left\{\ba{cl} x&\text{for }x \leq 0,\\ \frac{1{+}\eps}\eps \,x
  &\text{for } x \in[0,\eps],\\ x{+}1&\text{for }x\geq \eps; 
\ea\right. \quad \text{and}\quad 
X_\eps(y)=\left\{\ba{cl}  y&\text{for }y\leq 0, \\
  \frac{\eps}{1{+}\eps}\, y &\text{for }y \in [0,1{+}\eps], \\
y{-}1&\text{for }y\geq 1{+}\eps.\ea\right. 
\]
For $u:[0,T]\ti\Omega \to \R$ and  $y \in \wh\Omega:={]{-} 1,2[}$ we
define the functions  
\[
U_\eps(t,y)=u(t,X_\eps(y)),\quad W_\eps(y)=w_\eps(X_\eps(y)), \quad 
A_\eps(y)=\frac{a_\eps(X_\eps(y))}{X'_\eps(y)},
\]
and the functionals $\wh {\mathscr D}_\eps$ via
${\mathscr D}_\eps(u)=\wh{\mathscr D}_\eps(U_\eps)$ and find
\[
\wh{\mathscr D}_\eps(U):=\int_0^T \int_{\wh\Omega} \frac1{2 A_\eps U} \wh
I_\eps[\dot U]^2 + \frac{A_\eps U}{2} \Big( \big( \log(
U/W_\eps)\big)'\Big) ^2 \dd y \dd t,
\]
where $\wh I_\eps[\dot U](Y)= \int_{-1}^y \dot U(\eta) X'_\eps(\eta)
\dd \eta = - \int_y^2 \dot U(\eta) X'_\eps(\eta) \dd \eta$. 

Following the arguments in \cite[Sec.\,3.2]{Lier12VME} it is not
difficult to establish the $\Gamma$-convergence of $\wh{\mathscr D}_\eps$
to $\wh{\mathscr D}_0$, where the latter is given in the form 
\begin{align*}
&\wh{\mathscr D}_0(U):=\int_0^T \int_{\wh\Omega} \frac1{2 \wh A U} \wh
I_0[\dot U]^2 + \frac{\wh A  U}{2} \Big( \big( \log(
U/\wh W)\big)'\Big) ^2 \dd y \dd t\\
&\text{with }\ \wh I_0[\dot U](y):= \int_{-1}^y \dot U(\eta) \wh M(\eta)
\dd \eta = - \int_y^2 \dot U(\eta)\wh M(\eta) \dd \eta,  \\
&\text{where }(\wh A(y), \wh W(y),\wh M(y)):= \left\{\ba{cl}
  (a_-(y), \ee^{-V_-(y)}/Z_0, 1) &\text{for }y<0,\\ 
  (a_*(y), \ee^{-V_*(y)}/Z_0, 0) &\text{for }y\in[0,1],\\ 
  (a_+(y), \ee^{-V_+(y)}/Z_0, 1) &\text{for }y>1.
\ea \right.  
\end{align*}

\textit{Step 2. Minimization over the rescaled layer:}\/ 
The main structure in this limit model is that $\wh{\mathscr D}_0$ does not
depend on $\dot U(t,\cdot)|_{{]0,1[}}  $, since $\wh M(y)=0$ for $y \in
[0,1]$. Moreover, on this interval $\wh I_0[\dot U(t,\cdot)]$ is
constant, namely $\mu_u(t):= \int_{-1}^0 \dot U(\eta)\dd \eta= \wh
I_0[\dot U(t)](y)$ for  $y\in [0,1]$.
Thus, given the value $\mu_u(t)$ we can obtain the optimal profile of
$U(t)|_{[0,1]}$ from the boundary values $U(t,0)$ and $U(t,1)$ and
minimizing the functional $\wh\calG(\mu_u(t);\cdot)$ given via 
\begin{equation}
  \label{eq:calG}
  \wh\calG(\alpha, U):=\int_0^1 \frac{\alpha^2}{2\wh A U} + \frac{\wh A U}{2}
  \Big( \big( \log(U/\wh W)\big)'\Big)^2 \dd y.
\end{equation}
Now, Proposition \ref{pr:Layer} in Appendix \ref{App:A} provides the
explicit formula 
\begin{align*}
   \wh G(\alpha,u_0,u_1)&:=\min\set{ \wh\calG(\alpha, U) }{ U>0, \
    U(0)=u_0, \ U(1)=u_1}
  \\
  &= A_* \sqrt{\tfrac{u_0u_1}{w_-w_+}} \CCC\Big(\tfrac1{A_*}
  \sqrt{\tfrac{w_-w_+}{u_0u_1}}\,\alpha \Big) + A_*
  \sqrt{\tfrac{u_0u_1}{w_-w_+}} \CCC^*\Big( \log\big(\tfrac{ u_0
    w_+}{u_1 w_-}\big) \Big)
\end{align*}
with $w_-=\wh W(0)=w_0(0^-)$, $w_+=\wh W(1)=w_0(0^+)$, and $A_*=
\big(\int_0^1 1/(\wh A(y)\wh W(y)) \dd y\big)^{-1}$. Inserting the
definitions of $\wh A$ and $\wh W$ gives exactly the formula for $A_*$
in the theorem.

Thus, we have constructed a simpler functional $\ol{\mathscr D}_0$, which is
given by
\begin{align*}
\ol{\mathscr D}_0(U):=\int_0^T \bigg[\int_{{]-1,0[}\cup{]1,2[}} \Big(\frac1{2 \wh A U} \wh
I_0[\dot U]^2 + \frac{\wh A  U}{2} \Big( \big( \log(
U/\wh W)\big)'\Big) ^2  \Big)  \dd y &  \\[-0.5em]
+ \wh G\Big(\wh I_0[\dot U](0), U(t,0),U(t,1) \Big) & \bigg] \dd t  ,
\end{align*}
satisfies the lower bound $\wh{\mathscr D}_0(U)\geq \ol {\mathscr D}_0(U)$ for all
$U$, and has the important property that it does not depend on
$U|_{[0,T]\ti{]0,1[}}$. 

\textit{Step 3. Relation between $\ol{\mathscr D}_0$ and ${\mathscr D}_0$:}\/
Using the special form of $\calR_0$ and $\calR_0^*$ stated in
\eqref{eq:Memb.calR0} and the theorem, we define the limiting
dissipation functional ${\mathscr D}_0(u)= \int_0^T  \calR_0(u,\dot u ) +
\calR_0^*(u,{-}\rmD \calE_0(u) \big) \dd t$.
By construction and the special form of $\wh G$ the functional
$\ol{\mathscr D}_0$ is closely related to ${\mathscr D}_0$ in
the following way. For any function $U:[0,T]\ti [{-}1,2]\to \R$,
we may define $u:[0,T]\ti [{-}1,1]$ by $u(t,x)=U(t,Y_0(x))$,
where $Y_0(x)=x $ for $x\leq 0$ and $Y_0(x)=x{+}1$ otherwise. Then, we
have ${\mathscr D}_0(u)=\ol{\mathscr D}_0(U)\leq \ {\mathscr D}_0(U)$. 

Moreover, for any function $u$ one can construct an optimal 
$U$ as follows. We split $u$ at $x=0$, the right part is shifted by 1 to
the right, and the minimizer of $U(t,\cdot) \in \rmH^1([0,1])$ of
$\wh\calG( I[\dot u(t)](0),u(t,0^-),u(t,0^+))$ is inserted into the
gap. Then, ${\mathscr D}_0(u)=\wh {\mathscr D}_0(U)$.

\textit{Step 4. The liminf estimate \eqref{eq:defEDP.c}:}\/ 
To establish the fundamental liminf
estimate we consider, w.l.o.g., sequences $u_\eps$
satisfying $1/R \leq u_\eps \leq R$ for some large $R>1$.
In particular, by minimum and maximum principles these bounds
can be expected for solutions of \eqref{eq:Memb-1}.
Defining $U_\eps(t,y)=u_\eps(t,X_\eps(y))$ 
we again have $U_\eps(t,y)\in [1/R,R]$. Thus, we find a
subsequence $\eps_k\to 0$ such that
\[
u_\eps \weak u_0 \text{ in } \rmL^2([0,T]\ti \Omega) \quad\text{and}
\quad U_\eps \weak U_0 \text{ in } \rmL^2(0,T;\rmL^2(\wh\Omega)). 
\]
Moreover, we have $u_0(t,x)=U_0(t,Y_0(x))$. Now, using ${\mathscr
  D}_\eps(u_\eps)=\wh{\mathscr D}_\eps(U_\eps)$ we arrive at the
desired liminf estimate
\begin{equation}
  \label{eq:MembLiminf}
  \liminf_{\eps \to 0} {\mathscr D}_\eps(u_\eps) = \liminf_{\eps \to 0}
  {\mathscr D}_\eps(U_\eps) \geq \wh {\mathscr D}_0(U_0) \geq 
  \ol{\mathscr D}_0(U_0)= {\mathscr D}_0(u_0).  
\end{equation}
This concludes the proof of Theorem \ref{th:Membrane}. 
\hfill\rule{0.43em}{0.43em}\medskip

We conclude this section by observing that the EDP-limit of the
thin-layer diffusion system given by the classical GS
$(\PROB(\Omega),\calE_\eps,\calR_\eps^*)$ is a the generalized
GS for the membrane problem. For $\eps >0$ and for $\eps=0$ the
gradient structures are exactly the ones obtained from the
large-deviation principle, see Section \ref{sss:MP-Memb}. Hence, we
again found an instance where the diagram in Figure \ref{fig:LDPvsEDP}
commutes, that means that applying the large-deviation principle can
be interchanged with taking the EDP-limit $\eps \to 0$.
 
\section{From diffusion to reaction}
\label{se:D2R}

In our second major application of EDP-convergence as a microscopic
origin of generalized GS, we continue the work in
\cite{PeSaVe10DRGC,PeSaVe12CRGL,AMPSV12PLWG} which show that linear
reactions can be obtained as limits of diffusion for a suitably scaled
energy barrier.  In \cite{PeSaVe10DRGC,PeSaVe12CRGL} the method relies
on a quadratic energy functional and a classical gradient
structure. In \cite{AMPSV12PLWG} the pE-convergence for the entropic
GS is shown, but only diffusion along the reaction path is
allowed. In fact, the result therein gives EDP-convergence, if one 
takes the addition in \cite[Prop.\,4,4]{MiPeRe14RGFL} into account.

Here we generalize the latter work by also allowing diffusion in a
physical space $\Omega$, such that the resulting limit equation will
be a (linear) reaction-diffusion system.  Our physical domain
$\Omega\subset \R^d$ is bounded and has a Lipschitz boundary. For the
reaction path we choose $\Upsilon=[0,7] \subset \R$ and define the
cylinder $Q = \Omega \ti \Upsilon$.  (Indeed, $\Upsilon$ could by any
bounded or unbounded interval.)

\begin{figure}
\centering
\begin{tikzpicture}
\draw[->] (-0.5,0) -- (7.5,0) node[right] {$y$};
\draw[->] (0,-0.5) -- (0,3.2) node [right] {{\color{blue}$V(y)$}};
\node [below left] at (0,0) {$0\!$};
\node [below] at (2,0) {$2$};
\node [below] at (5,0) {$5$};
\node [below] at (6,0) {$6$};
\node [below right] at (7,0) {$\!7$};
\foreach \x in {1,...,6} \draw (\x, 0.1) -- (\x, -0.1); 
\draw[|-|,line width=0.2em] (-0.02,0)--(7,0) node[pos=0.6,above] {$\Upsilon$};

\draw[line width=0.4em,draw=blue] (0,1.5) .. controls (1.5,0) .. (2,0); 
\draw[line width=0.4em,draw=blue] (2,0) .. controls (2.9,0) 
  and (4,3) .. (5,3) ; 
\draw[line width=0.4em,draw=blue] (5,3) .. controls (5.7,3) 
  and (5.3,0) .. (6,0); 
\draw[line width=0.4em,draw=blue] (6,0) .. controls (6.3,0) .. (7,1); 
\end{tikzpicture}
\caption{The potential $V$ along the reaction path $\Upsilon =[0,7]$.}
\label{fig:Potential}
\end{figure}
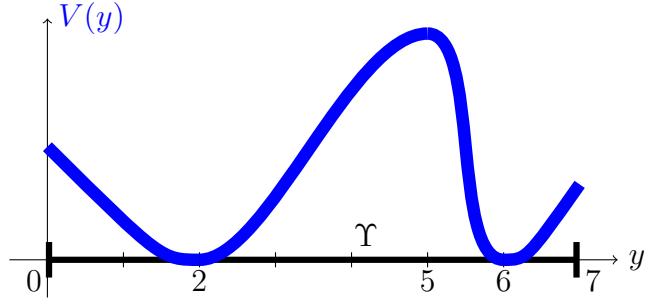

For densities $u\in \rmL^1(Q)$ the integral $\int_D \int_{y_0}^{y_1} u
\dd y \dd x $ denotes the number of particles per unit volume that are
in the subdomain $D \subset \Omega$ and have a reaction state $y \in
[y_0,y_1]\subset \Upsilon$. The evolution of the density $u$ is driven
by diffusion in the $x$-direction with diffusion constant $m_\Omega>0$
and a much faster diffusion in the $y$-direction with diffusion
constant $\tau_\eps \gg 1 $ to allow the particles to overcome a huge
potential barrier given by $V_\eps(y)=\frac1\eps V(y)$, see Figure
\ref{fig:Potential}.

 For simplicity we assume that the total
mass $\int_Q u\dd x \dd y$ as well as the volume $|\Omega|$ of the
physical domain equal 1. Hence, we can again consider the model as a
Markov process with continuous 
paths $t\mapsto (X_t,Y_t) \in \Omega\ti \Upsilon=Q$, whose distribution laws
can be described by densities $u(t)\in \PROB(Q)$.
The Kolmogorov forward equation  reads
  \begin{equation}
  \label{eq:D2R.1}
  \dot u =m_\Omega \Delta_x u + \tau_\eps\pl_y \big( \pl_y u +
  u\,\pl_yV_\eps \big),\quad 
  (\nabla_x u,\pl_y u {+} u \,\pl_yV_\eps )\cdot \nu=0 \text{ on }\pl Q.
\end{equation}
Clearly, the unique steady state $\wt w_\eps$ is independent of $x$
and takes the form
\[
\wt w_\eps(x,y)=w_\eps(y):= \frac1{z_\eps} \exp\Big({-}\frac1\eps V(y)\Big)
 \text{ with }  z_\eps:=\int_\Upsilon\exp\Big({-}\frac1\eps V(y)\Big)\dd y. 
\]
Equation \eqref{eq:D2R.1} is a Fokker-Planck equation and, hence, has the
Wasserstein gradient structure introduced in \cite{JoKiOt98VFFP} with
\begin{equation}
  \label{eq:D2R.F.R}
  \calE_\eps(u)=\iint_Q \lambda_\rmB\big(\frac{u}{w_\eps}\big)
 w_\eps \dd y \dd x \  \text{ and } \ 
\calR_\eps^*(u,\xi) = \iint_Q \Big(\frac {m_\Omega}2|\nabla_{\!x}\xi|^2 + 
 \frac{\tau_\eps} 2 (\pl_y\xi)^2 \Big)u \dd y \dd x. 
\end{equation}

For studying the limit $\eps\to 0$ we now assume that $V\in
\rmC^2(\Upsilon)$ has exactly two non-degenerate minimizers 
as pure states, where $V=0$ w.l.o.g, and one global maximum 
as barrier, namely
\begin{align} \label{eq:D2R.1c}
 &V(2)=V(6)=0, \ V(y)>0 \text{ on } \Upsilon\setminus\{2,6\}, \ 
   V''(2)>0,\ V''(6)>0; \\
 &V(5)> V(\ol y) \text{ on } \Upsilon\setminus\{5\}, \ V''(5)<0, 
\end{align}
see Figure \ref{fig:Potential}. 
(Again, any two points in $\Upsilon$ could be taken as the pure
states, and any point in between as barrier.)  As a consequence
$w_\eps$ concentrates in the points $y=2$ and $y=6$ in the limit, 
viz.
\begin{equation}
  \label{eq:w-eps20}
w_\eps \weaks w_0 = \alpha_0\delta_2 + \alpha_1\delta_6\in 
 \PROB(\Upsilon) , \ \ 
 \alpha_0 = \tfrac{\sqrt{V''(6)}}{\sqrt{V''(2)}{+}\sqrt{V''(6)}},\ 
 \alpha_1 = \tfrac{\sqrt{V''(2)}}{\sqrt{V''(2)}{+}\sqrt{V''(6)}}.
\end{equation}
Here the convergence means $\int_\Upsilon \phi(y)w_\eps(y)\dd y \to
\alpha_0\phi(2) + \alpha_1 \phi(6)$ for all $\phi\in \rmC^0(\Upsilon)$.

The important point is now to choose the diffusion constant
$\tau_\eps$ sufficiently large such that the transitions between $y=2$
and $y=6$ can occur on times of order 1. According to Kramer's
rule (see e.g.\ \cite{AMPSV12PLWG}), this is achieved by choosing
$m_\Upsilon>0$ and setting
\[
\tau_\eps := m_\Upsilon \int_\Upsilon \frac1{w_\eps(y)} \dd y, \text{
  with gives }  \frac{\tau_\eps} \eps\exp\big({-}V(5)/\eps \big)  \to m_\Upsilon 
  \tfrac{2\pi\,\big(\sqrt{V''(2)}+\sqrt{V''(6)}\big)}
   {\sqrt{-V''(5)}\:\sqrt{V''(2)V''(6)}}\ > 0.
\]

From the concentration of $w_\eps$ in the points $\{2,6\}$ we obtain
that $\wt w_\eps\in \PROB(Q)$ concentrates in the sets
$\Omega\ti\{2\}$ and $\Omega\ti \{6\}$, namely
\[
\wt w_\eps \ \weaks \ \wt w_0 := \chi_\Omega \otimes w_0 \text{ in } \PROB(Q).
\]
Recalling that $\calE_\eps$ is the relative entropy with respect to
$w_\eps$, the $\Gamma$-convergence  $\calE_\eps \Gweaks \calE_0$ appears
natural. To be more precise concerning densities and measures, 
we define 
\begin{align} \nonumber
\calE_0(\mu) &:= \left\{ \ba{cl}  \bfE((c_0,c_1)) &\text{for } 
        \mu=c_0\dd x{\otimes } \delta_2 +c_1\dd x{\otimes } \delta_6, \\ 
      \infty &\text{otherwise},  \ea \right.
\\  \label{eq:D2R.F0}
\text{where }&\bfE((c_0,c_1)):= \int_\Omega \Big( 
        \lambda_\rmB(\tfrac{c_0(x)}{\alpha_0}\big)\alpha_0  +  
        \lambda_\rmB(\tfrac{c_1(x)}{\alpha_1}\big)\alpha_1\Big) \dd x.   
\end{align}

\begin{proposition}\label{pr:D2R-calF}
We have $\calE_\eps \Gweaks \calE_0$ in the weak$^*$ topology of $\PROB(Q)$.
\end{proposition}
\begin{proof}
The liminf estimate is established in
\cite[Lem.\,9.4.3]{AmGiSa05GFMS}.  

To construct recovery sequences, we may restrict to the case
$\calE_0(\mu)<\infty$, since otherwise the liminf estimate provides
the result. Hence, we may assume $ \mu=c_0\dd x{\otimes } \delta_2
+c_1\dd x{\otimes } \delta_6$ and, using a nonnegative, continuous cut-off
function $\chi:y\mapsto \max\{1{-}|y|,0\}$, we can  define the measures 
\[
\mu_\eps =u_\eps(x,y) \dd x{\otimes}\dd y \text{ with }
 u_\eps(x,y)=
   c_0(x)w_\eps(y)\beta_{0,\eps}\chi(y{-}2) +
c_1(x) w_\eps(y) \beta_{1,\eps}\chi(y{-}6),
\]
where the normalization constants $\beta_{j,\eps}$ are given by 
\[ \ts
\beta_{0,\eps}\int_\Upsilon w_\eps(y)\chi(y{-}2)\dd y =1 
=\beta_{1,\eps}\int_\Upsilon w_\eps(y)\chi(y{-}6)\dd y,
\]
which implies $\beta_{j,\eps}\to 1/\alpha_j>0$ for $\eps \to 0$. 
 Then, we easily find $\calE_\eps (u_\eps)\to \calE_0(\mu)$. 
\end{proof}

Thus, the limit evolution will be described by the densities $c_0$ and
$c_1$ on $\Omega$ for the particles being in the pure states $y=2$ and
$y=6$, respectively. In particular, in the limit $\eps\to 0$ the time
that the particles spend along the reaction path away from these
points, i.e.\ in $\Upsilon\setminus\{2,6\}$, is $0$.
   
One difficulty in deriving the liminf estimate for 
De Giorgi's dissipation functional
\[
{\mathscr D}_\eps(u):= \int_0^T \calR_\eps(u(t),\dot u(t)) +
    \calR_\eps^*\big(u(t),{-}\rmD \calE_\eps(u(t))\big) \dd t
\]
is that $\calR_\eps$ is only implicitly defined
via the Legendre transform of $\calR_\eps^*$. Moreover, we are not able
to employ the classical Wasserstein gradient flow theory in
\cite{AmGiSa05GFMS} using the Benamou-Brenier formulation, because of
the different roles of the diffusion in $x$ with mobility $m_\Omega$
and the diffusion in $y$ with mobility $\tau_\eps \to \infty$. The
first step to establish the following result follows the idea in
\cite{MaaMie15?GSRC1}, where one obtains a lower estimate by replacing
$\calR_\eps(u,\dot u)$ by the smaller term $\langle \xi_\eps,\dot
u_\eps\rangle - \calR_\eps^*(u_\eps,\xi_\eps)$ and by choosing a
suitable recovery sequence $\xi_\eps \to \xi_0$ for the limit passage
$\eps \to 0$. Finally, one takes the supremum over all $\xi_0$ to
recover $\calR_0$ as dual of $\calR_0^*$. The second step involves a
suitable transformation of the reaction variable $z=Z_\eps(y)$ (first
introduced in \cite{AMPSV12PLWG}) which
allows us to control the relative densities $v_\eps:=u_\eps/w_\eps$
and the dual potentials $\xi_\eps$ along the reaction path $\Upsilon$.

In the following result we will again describe the limit GS
$(\PROB(Q),\calE_0,\calR_0)$ by a reduced GS
$(\PROB(\Omega\ti\{0,1\}), \bfE,\bfR)$, since in the limit every 
$\mu \in \PROB(Q)$ with finite relative entropy
satisfies $\mu = c_0\dd x {\otimes}\delta_2(y) +
c_1\dd x {\otimes}\delta_6(y) $ with $(c_0\dd x,c_1\dd x) \in
\PROB(\Omega\ti\{0,1\})$, see \eqref{eq:D2R.F0}.

\begin{theorem}[From diffusion to reaction-diffusion]\label{th:D2R-calJ}
  The family of gradient systems $(\PROB(Q),\calE_\eps,\calR_\eps)$
  defined via \eqref{eq:D2R.F.R} converges in the EDP sense to the
  gradient system $(\PROB(Q),\calE_0,\calR_0)$, where $\calE_0$ is
  given in \eqref{eq:D2R.F0} via $\bfE$ and accordingly $\calR_0$ is
  given via $\bfR$, which is defined in terms of the dual dissipation
  potential
\[
\bfR^*(\bfc,\bfeta):= \int_\Omega \frac{m_\Omega}2 \big(c_0|\nabla_{\!x}
\eta_0|^2 {+} c_1|\nabla_{\!x}\eta_1|^2\big) +  m_\Upsilon 
\sqrt{\tfrac{c_0c_1}{\alpha_0\alpha_1}} \CCC^*(\eta_1{-}\eta_0) 
\dd x.
\]
\end{theorem}

The above result means that the limiting GS is a generalized gradient
system defined for $\bfc=(c_0,c_1)\in \PROB(\Omega{\ti}\{0,1\})$,
where the limiting system is the coupled system of linear PDEs given
in the form
\[
\dot c_0= m_\Omega \Delta c_0 - m_\Upsilon 
         \big( c_0/\alpha_0 - c_1/\alpha_1\big),  \quad
\dot c_1= m_\Omega \Delta c_1 + m_\Upsilon 
         \big( c_0/\alpha_0 - c_1/\alpha_1\big),
\]
with Neumann boundary conditions $\nabla c_j \cdot \nu=0$.  We
emphasize that the original GS
$(\PROB(Q),\calE_\eps,\calR_\eps)$ is the classical GS for
the Fokker-Planck equation, while the  EDP limit provides the 
generalized gradient structure discussed in Section \ref{sss:LinRDS}. 
We observe that for $\eps>0$ as well as for $\eps=0$ we have the GS
that is induced by the large-deviation principle discussed in Section
\ref{ss:Markov}.  Thus, we have found another instance of the
interchangeability of the large-deviation principle and the EDP-limit,
as displayed in Figure \ref{fig:LDPvsEDP}.\bigskip 

\noindent\emph{Sketch of proof of Theorem \ref{th:D2R-calJ}:} 
Since the $\Gamma$-convergence $\calE_\eps \Gweaks \calE_0$ was
already established in Proposition \ref{pr:D2R-calF}, it remains to
show the liminf estimate for the dissipation functional ${\mathscr D}_\eps$.
More precisely, assume $u_\eps(t)\weaks \mu(t)=c_0(t)\dd
x{\otimes}\delta_2(y) + c_1(t)\dd x {\otimes} \delta_6(y)$ in
$\PROB(Q)$ for all $t\in [0,T]$ such that $\sup_{t\in [0,T]}
\calE_\eps(u_\eps) < \infty$; then, we have to show
\begin{equation}
  \label{eq:bfJ(bfc)}
\liminf_{\eps\to 0} {\mathscr D}_\eps(u_\eps) \geq {\bfD}(\bfc):= \int_0^T
\bfR(\bfc,\dot \bfc) + \bfR^*(\bfc,{-}\rmD\bfE(\bfc))\dd t.
\vspace*{0.3em}
\end{equation}

\textit{Step 1. Dualization of $\calR_\eps$:}\/ 
The first major idea follows \cite{MaaMie15?GSRC1} and exploits the definition
of $\calR_\eps$ as Legendre transform of 
$\calR_\eps^*$. Introducing the functional
\[
{\mathscr B}_\eps(u,\xi):= \int_0^T \langle \xi,\dot u\rangle - \calR^*(u,\xi)
 + \calR_\eps^*(u,{-}\rmD\calE_\eps(u))\dd t,
\]
we easily see that ${\mathscr D}_\eps(u)$ can be reconstructed via $\sup_\xi
{\mathscr B}_\eps(u,\xi)$.  Using the definitions of $\calE_\eps$ and
$\calR_\eps^*$ we have the explicit form
\begin{align*}
{\mathscr B}_\eps(u,\xi)=\int_0^T\hspace{-0.5em} \int_Q &\Big[\xi\dot u - 
 \frac{m_\Omega}2|\nabla_{\!x}\xi|^2u - \frac{\tau_\eps}2 (\pl_y\xi)^2 u
 \\&\quad+\frac{m_\Omega}2\,\frac{|\nabla_{\!x}u|^2}u  +
  \frac{\tau_\eps}2\Big( \pl_y\big(\log(u/w_\eps)\big)\Big)^2 u\Big] 
  \dd y \dd x \dd t.
\vspace{0.3em} 
\end{align*}

\textit{Step 2. Rescaling the reaction-path variable.}\/ 
The second major idea follows \cite[Sec.\,2.1]{AMPSV12PLWG}, 
where no $x$-direction was present. We
define the diffeomorphism $Z_\eps :\Upsilon\to\calZ:=[0,1]$ and its
inverse $Y_\eps=Z_\eps^{-1}:\calZ\to \Upsilon$ via
\[
z=Z_\eps(y):=\frac{m_\Upsilon}{\tau_\eps} \int_{\ol y=0}^y \frac1{w_\eps(\ol y)} 
   \dd \ol y \ \text{ and } \ 
Y'_\eps(z)=\frac{\tau_\eps}{m_\Upsilon} w_\eps( Y_\eps(z)). 
\]
The transformed equilibrium density $\wh w_\eps$ on $\calZ$ is 
\begin{equation}
  \label{eq:wh-w-eps20} 
\wh w_\eps(z) :=w_\eps(Y_\eps(z))Y'_\eps(z) \ \text{ and satisfies }
\wh w_\eps \,\weaks \, \wh w_0 := \alpha_0 \delta_0 + \alpha_1 \delta_1. 
\end{equation}
Indeed, for the latter statement we first use that for all $g\in
\rmC^0(\calZ)$ we have the identity $\int_\calZ g(z)\wh w_\eps(z)\dd z =
\int_\Upsilon g(Z_\eps(y)) w_\eps(y) \dd y$.  Recalling that $V$
has a unique global maximum at $y=5$, the function $Z_\eps$ converges
uniformly on compact subsets of $\Upsilon\setminus\{5\}$ to the step
function $Z_0(y)=0$ for $y<5$ and $Z_0(y)=1$ for $y>5$. With this and 
\eqref{eq:w-eps20} we
conclude $ \int_\Upsilon g(Z_\eps(y)) w_\eps(y) \dd y \to \alpha_0
g(Z_0(2)) + \alpha_1 g(Z_0(6))$ which is the desired result
\eqref{eq:wh-w-eps20}.

To estimate ${\mathscr B}_\eps$ in the limit $\eps\to 0$ we use now the
independent variable $z=Z_\eps(y)$ and the dependent variables
\[
v(t,x,z)=\frac{u(t,x,Y_\eps(z))}{w_\eps(Y_\eps(z))} \ \text{ and } \ 
\zeta(t,x,z)=\xi(t,x,Y_\eps(z)).
\]
Introducing the domain $\wh Q=\Omega \ti \calZ$ we find ${\mathscr B}_\eps(u,\xi)=\wh {\mathscr B}_\eps(v,\zeta)$ with 
\[
\wh{\mathscr B}_\eps(v,\zeta)= \int_0^T\hspace*{-0.5em}\int_{\wh Q} 
  \zeta \dot v\wh w_\eps 
  - \frac{m_\Omega}2|\nabla_{\!x}\zeta|^2 v \wh w_\eps
  - \frac{m_\Upsilon}2 (\pl_z\zeta)^2 v
  + \frac{m_\Omega}2\,\frac{|\nabla_{\!x} v|^2}{v}  \wh w_\eps
  + \frac{m_\Upsilon}2\,\frac{ (\pl_z v)^2}{v} \;\dd z \dd x \dd t.
\]
The transformation of ${\mathscr B}_\eps$ to $\wh{\mathscr B}_\eps$ follows easily 
by using the relations 
\[
u \dd y = v \wh w_\eps \dd z, \quad \pl_y\xi = \frac{\pl_z \zeta}{Y'_\eps(z)}, 
\quad \text{and }
\frac{\tau_\eps}{(Y'_\eps(z))^2} = \frac{m_\Upsilon}{\wh w_\eps(z)}.
\vspace{0.3em} 
\]

\textit{Step 3. The $\Gamma$-limit for $\wh {\mathscr B}_\eps(\cdot,\zeta)$:}\/
The importance of the new form $\wh{\mathscr B}_\eps$ is that the dependence
on $\eps$ only occurs in the weighting measure $\wh w_\eps$. Since
$\wh w_\eps$ concentrates in the points $z=0$ and $1$, the three terms
that are multiplied by the weight $\wh w_\eps$ will converge 
to simple integrals over $[0,T]\ti\Omega$ for the densities
$c_0$ and $c_1$ respectively. In contrast there are two terms not
involving $\wh w_\eps$, but these terms only involve derivatives in
the $z$-direction. In particular, they control the smoothness of
$\zeta$ and $v$ in $z$-direction, namely $\sqrt{v} \in
\rmL^2([0,T]\ti\Omega; \rmH^1(\calZ))$ such that $v\wh w_\eps$ indeed
has a well-defined limit. With this and $\wh w_\eps \weaks \wh
w_0=\alpha_0\delta_1 + \alpha_1\delta_1$ (cf.\ \eqref{eq:wh-w-eps20}),  
it is possible to show that for fixed and sufficiently smooth $\zeta$
we have  $\wh{\mathscr B}_\eps(\cdot,\zeta)
\Gweak \wh{\mathscr B}_0(\cdot,\zeta)$ with 
\begin{align*}
\wh{\mathscr B}_0(v,\zeta)= \int_0^T\hspace*{-0.5em}\int_\Omega \bigg[&\int_\calZ 
\frac{m_\Upsilon}2\Big(\frac{ (\pl_z v)^2}{v} - (\pl_z\zeta)^2
v\Big)\dd z  \\
& + \sum_{j=0}^1  \alpha_j \Big(\zeta_j\dot v_j 
        -\frac{m_\Omega}2 v_j|\nabla_{\!x} \zeta_j|^2 
        + \frac{m_\Omega}2\,\frac{|\nabla_{\!x} v_j|^2}{v_j} 
      \big)  \:\bigg]\:\dd x \dd t,
\end{align*}
where $v_j(t,x)=v(t,x,j)$ and $\zeta_j(t,x)=\zeta(t,x,j)$ for $j=0$
and $j=1$.\vspace{0.3em} 

\textit{Step 4. Minimization over the reduction path profile.}\/ 
Note that in the definition of $\wh{\mathscr B}_0$, the values of the 
functions $v$ and $\zeta$ for $z\in {]0,1[}$ only occur in the first
integrand (with factor $m_\Upsilon$). Hence, one can eliminate the
integral by  taking the supremum in $\zeta$ and the infimum in $v$ 
for given  boundary values at $z=0$ and $z=1$. The relevant functional
reads 
\[
\calN(v,\zeta)=\int_0^1 \Big(\frac{v'(z)^2}{2v(z)} - \frac 12 \zeta'(z)^2
v(z) \Big) \dd z,
\]
and Proposition \ref{pr:calN} provides the following explicit inf-sup formula
\begin{align*}
&\inf\Bigset{ \sup\bigset{\calN(v,\zeta)}{\zeta(0)=\zeta_0,\:\zeta(1)=\zeta_1}
}{ v(0)=v_0,\:v(1)=v_1,\: v>0 } \ \\
&= \ \sqrt{v_0v_1}\,\CCC^*(\log
v_1{-}\log v_0) - \sqrt{v_0v_1}\, \CCC^*(\zeta_1{-}\zeta_0)=: 
N(\zeta_1{-}\zeta_0, v_0,v_1). 
\end{align*}
Thus, we can reduce $\wh {\mathscr B}_0$ to a functional $\bfB$ on
 $\bfv=(v_0,v_1)$ and $\bfzeta= (\zeta_0,\zeta_1)$, namely  
\begin{align*}
\bfB(\bfv,\bfzeta):= \int_0^T\hspace*{-0.5em}\int_\Omega \bigg[& m_\Upsilon  
           N(\zeta_1{-}\zeta_0, v_0,v_1)   \\
&      +\sum_{j=0}^1 \alpha_j \Big(\zeta_j\dot v_j 
        -\frac{m_\Omega}2 v_j|\nabla_{\!x} \zeta_j|^2 
        + \frac{m_\Omega}2\,\frac{|\nabla_{\!x} v_j|^2}{v_j} 
      \Big)  \bigg]\dd x \dd t.
\end{align*}
The inf-sup definition of $N$ provides the following relation between
$\wh{\mathscr B}_0$ and $\bfB$:
\begin{align}
  \label{eq:calB.bfB}
&\forall\ v\text{ with }v|_{[0,T]\ti\Omega\ti\{0,1\}}=\bfv\ 
\exists\ \zeta\text{ with }\zeta|_{[0,T]\ti\Omega\ti\{0,1\}}=\bfzeta:
\quad \wh{\mathscr B}_0(v,\zeta)\geq \bfB(\bfv,\bfzeta). 
\vspace{0.3em}
\end{align}

\textit{Step 5. Identification of the limits:}\/ 
It now remains to relate the limit functions $\bfv$ to the weak limit
of the sequence $u_\eps$. For this, we consider a sequence $u_\eps$ as
in \eqref{eq:defEDP.c}, i.e.\ $u_\eps \weaks u$ and
$\calE_\eps(u_\eps(t)) \leq C<\infty$. By the definition of $v$ we have
\begin{equation}
  \label{eq:D2R-10}
  \int_0^T\hspace{-0.5em}\int_Q u_\eps(t,x,y)\phi(t,x,y)\dd y\dd x \dd t = 
\int_0^T\hspace{-0.5em}\int_Q v_\eps(t,x,Z_\eps(y)) w_\eps(y)
\phi(t,x,y)\dd y \dd x \dd t. 
\end{equation}
Without loss of generality we assume $\infty > C\geq
{\mathscr D}_\eps(u_\eps)=\sup_\xi {\mathscr B}_\eps(u_\eps, 
\xi) = \sup_\zeta \wh {\mathscr B}(v_\eps,\zeta)$, which gives the bound 
\begin{equation}
  \label{eq:sqrt-v}
\| \sqrt{v_\eps} \|_{\rmL^2([0,T]\ti\Omega;\rmH^1(\calZ))} \leq C.
\end{equation}
This implies H\"older continuity of $v(t,x,\cdot):\calZ\to \R$.  Moreover,
$Z_\eps$ converges uniformly to $0$ and $1$ near $y=2$ and
$y=6$.  Hence, we can pass to the limit in
\eqref{eq:D2R-10} and obtain
\begin{align*}
&\int_0^T\hspace*{-0.5em}
\int_\Omega c_1(t,x)\phi(t,x,2){+}c_2(t,x)\phi(t,x,6)\dd x \dd t 
=\lim_{\eps\to 0}  \int_0^T\hspace*{-0.5em}
\int_Q u_\eps(t,x,y)\phi(t,x,y)\dd y\dd x \dd t \\
&=\lim_{\eps\to 0}  \int_0^T\hspace*{-0.5em}
\int_Q v_\eps(t,x,Z_\eps(y) w_\eps(y)\phi(t,x,y)\dd y\dd x \dd t 
= \int_0^T\hspace*{-0.5em}
\int_\Omega \sum\nolimits_{0}^1 v_j(t,x) \alpha_j \phi(t,x,2{+}4j) \dd x \dd t,    
\end{align*}
which means $\bfc=(c_0,c_1)= (\alpha_0 v_0, \alpha_1v_1)$. Using the
explicit form of $\bfR^*$ and $\bfE$ implies
\begin{equation}
  \label{eq:bfB-bfR}
  \int_0^T \langle \bfzeta,\dot\bfc\rangle {-} \bfR^*(\bfc,\bfzeta){+}
\bfR^*(\bfc,-\rmD\bfE(\bfc)) \dd t =\bfB
\big((c_1/\alpha_0,c_2/\alpha_1), \bfzeta\big).
\vspace{0.3em}
\end{equation}

\textit{Step 6. The liminf estimate:}\/ With these preparations we can
now complete the liminf estimate. By the
construction of $v_\eps$ and $\xi_\eps(t,x,y)=\zeta(t,x,Z_\eps(y))$
we obtain the relations
\[
{\mathscr D}_\eps(u_\eps) \geq {\mathscr B}_\eps(u_\eps,\xi_\eps) =  \wh{\mathscr B}_\eps(v_\eps,\zeta),
\]
where $\zeta$ is now fixed.  For the sequence $u_\eps$ as given in
Step 5, we can further assume that $v_\eps \weak v$ in
$\rmL^2([0,T]{\ti}\Omega;\rmC^0(\calZ))$, using \eqref{eq:sqrt-v}.
According to Step 3 the liminf for $\eps \to 0$ yields
\[
\liminf_{\eps\to 0} {\mathscr D}_\eps(u_\eps) \ \geq \
\liminf_{\eps\to 0}\wh{\mathscr B}_\eps(v_\eps,\zeta) \ \geq \
\wh{\mathscr B}_0(v,\zeta) \ \geq \ 
\bfB\big((c_0/\alpha_0,c_1/\alpha_1),\bfzeta),
\]
where for the last estimate we have to choose $\zeta$ according to
\eqref{eq:calB.bfB} to fit the limit $v$ and $\bfzeta=
\zeta|_{[0,T]\ti\Omega\ti\{0,1\}}$. Nevertheless, the functions
$\bfzeta=(\zeta_0,\zeta_1)$ are still free.  Using the
characterization \eqref{eq:bfB-bfR} and taking the supremum over all
$\bfzeta$ gives the desired lower bound:
\begin{align*}
\liminf_{\eps\to 0} {\mathscr D}_\eps(u_\eps)\ &\geq \ \sup_\bfzeta
\wh{\mathscr B}_0(v,\zeta)  
\ \geq \ \sup_\bfzeta \int_0^T \!\!\big(\langle \bfzeta,\dot\bfc\rangle {-}
\bfR^*(\bfc,\bfzeta){+} 
\bfR^*(\bfc,-\rmD\bfE(\bfc))  \big)\dd t\\
&=\ \int_0^T\!\! \big(\bfR(\bfc,\dot\bfc) + \bfR^*(\bfc,-\rmD\bfE(\bfc))
\big)  \dd t\  = \ \ {\bfD}(\bfc). 
\end{align*}
Thus the desired estimate \eqref{eq:bfJ(bfc)} is established, which
finishes the proof of Theorem \ref{th:D2R-calJ}.\hfill\rule{0.43em}{0.43em}

\appendix
\section{Evaluation of some functionals}
\label{Appendix}

Here we give explicit calculations for the functional $\calG$ occurring
in the membrane limit and the functional $\calN$ occurring in the limit
of diffusion to reaction. It is surprising that both functional are
closely related, see \eqref{eq:calM.calG}. 

\subsection{Derivation of the potential $G(\alpha,u_0,u_1)$}
\label{App:A}

We first give the result of the standard case of constant
coefficients $\wh A$ and $\wh W$, which was already derived in
\cite[Prop.\,4.4]{MiPeRe14RGFL}  
under the restriction $u_0+u_1=1$. For  the functional 
\[
  \calG(\alpha, u):= \int_0^1 \frac{\alpha^2+ u'(x)^2}{2
  u(x)} \dd x,
\]
we define the value function 
\begin{equation}
  \label{eq:Galpha}
  G(\alpha,u_0,u_1):= \min\Big\{\; \calG(\alpha, u) \;\Big| 
\; u\in\rmH^1(0,1),\ u(0)=u_0,\ u(1)=u_1, \ u>0\;\Big\}, 
\end{equation}
and give a full proof of the
derivation of the explicit formula.

\begin{proposition}\label{pr:MMM}
For all $\alpha\in \R$ and $u_0,u_1$ we have 
\begin{equation}
  \label{eq:A.MMM}
\begin{aligned}
G(\alpha,u_0,u_1)&=\sqrt{u_0u_1}\, \CCC\Big(
\frac{\alpha}{\sqrt{u_0u_1}}\Big) +\sqrt{u_0u_1}\,\CCC^*(\log u_1{-}\log u_0), 
\end{aligned} 
\end{equation}
where the last term simplifies to $G(0,u_0,u_1)=2(\sqrt{u_0}{-}\sqrt{u_1})^2$.
Moreover, the unique minimizer is given by 
\begin{equation}
  \label{eq:A.Umin}
  u(x) = (1{-}x)u_0 + x u_1 + b(x^2{-}x)\quad \text{with } b=u_0+u_1 -
\sqrt{\alpha^2 {+} 4u_0u_1}. 
\end{equation}
\end{proposition}
\begin{proof} Since the integrand is convex, there is a unique
  minimizer $u$. Denoting the integrand by $f(u,u')$ the
  Euler-Lagrange equations $- \big(\pl_{u'}f(u,u')\big)'+ f(u,u')=0$
  are $u u'' - (u')^2 +\alpha^2=0$. By Noether's theorem we also have
  the first integral $u'\pl_{u'} f(u,u')-f(u,u')= \big((u')^2 {-}
  \alpha^2 \big)/(2u)=\gamma/2=$const. From $(u')^2=\alpha^2 + \gamma
  u$ it is now easy 
  to see that all solutions of the Euler-Lagrange equations are
  parabolas. Using the boundary conditions we find $u$ in
  \eqref{eq:A.Umin}, where $\gamma = 4b $.  

To evaluate the integral we restrict to the 
case $u'(x)>0$ on $[0,1]$, which means $2u_0 <
\sqrt{\alpha^2{+}4u_0u_1} < 2 u_1$. In the other cases, one can do the
calculation on all monotone parts in a similar fashion and add the
result. We use \eqref{eq:A.Umin} and $(u')^2=\alpha^2 +\gamma u$ to obtain 
\begin{align*}
\int_0^1 \frac{\alpha^2+ u'(x)^2}{2 u(x)} \dd x &
 = \frac{\gamma}2 + \alpha^2 \int_0^1 \frac{\rmd x} {u(x)}
  =2b + \alpha^2 \int_{u_0}^{u_1} \frac{\rmd u } {u\sqrt{\alpha^2{+} \gamma u}},  
\\
&=2b -2\alpha \mathop{\mathrm{\mathrm{artanh}}}\sqrt{1{+}\tfrac{\gamma}{\alpha^2} u_1}
   + 2\alpha \mathop{\mathrm{\mathrm{artanh}}}\sqrt{1{+}\tfrac{\gamma}{\alpha^2} u_0}.
\end{align*}
To proceed we first observe $\alpha^2 {+} \gamma u_j=
\big(\sqrt{\alpha^2{+}4u_0u_1} -2 u_j\big)^2$, which gives
$\sqrt{\alpha^2 {+} \gamma u_0} = \sqrt{\alpha^2{+}4u_0u_1} -2 u_0$
and $\sqrt{\alpha^2 {+} \gamma u_1} = 2u_1 - \sqrt{ \alpha^2{+}
  4u_0u_1}$. Now employing the addition rule $\mathop{\mathrm{\mathrm{artanh}}}
(x)+ \mathop{\mathrm{\mathrm{artanh}}} (y) =\mathop{\mathrm{arsinh}} \big((x{+}y)
/ \sqrt{(1{-}x^2) (1{-}y^2)}\big)$ and $\gamma=4b $ gives
\eqref{eq:A.MMM}.
\end{proof} 

In Section \ref{se:Membrane} we need a more general version with 
non-constant functions $\wh A$ and $\wh W$:
\begin{equation}
  \label{eq:calG.A}
  \wh\calG(\alpha, U):=\int_0^1 \frac{\alpha^2}{2\wh A U} + \frac{\wh A U}{2}
  \Big( \big( \log(U/\wh W)\big)'\Big)^2 \dd y.
\end{equation}
We will show that the influence of the coefficient functions $\wh A$
and $\wh W$ can be calculated from Proposition \ref{pr:MMM} by a
suitable rescaling of the layer variable in the form $x=X(y)$.

\begin{proposition}\label{pr:Layer}
We have the following formula:
\begin{align*}
  & \wh G(\alpha,u_0,u_1):=\min\set{ \wh\calG(\alpha, U) }{ U>0, \
    U(0)=u_0, \ U(1)=u_1}
  \\
  &= A_* \sqrt{\tfrac{u_0u_1}{w_0w_1}} \,\CCC\Big(\tfrac1{A_*}
  \sqrt{\tfrac{w_0w_1}{u_0u_1}}\,\alpha \Big) + A_*
  \sqrt{\tfrac{u_0u_1}{w_0w_1}} \,\CCC^*\Big( \log\big(\tfrac{ u_0
    w_1}{u_1 w_0}\big) \Big),
\end{align*}
where $w_0=\wh W(0)$, $w_1= \wh W(1)$, and $A_*=\mathrm{Harm}(\wh
A \wh W)=\big(\int_0^1 1/(\wh A(y) \wh W(y)) \dd y\big)^{-1}$.
\end{proposition}
\begin{proof} We define the new independent variable $z$ and a
new function $v(z)$ via 
\[
z=Z(y):=A_*\int_0^y \frac{\rmd \eta}{\wh A(\eta) \wh W(\eta)} \quad
\text{and} \quad v(Z(y))= \frac{U(y)}{ \wh W(y)},
\]
where $A_*=\big(\int_0^11/(\wh A(\eta) \wh W(\eta))\dd\eta\big)^{-1}
$.  By definition we have $Z(0)=0$ and $Z(1)=1$, and the inverse $Y$
of $Z$ maps $[0,1]$ into itself again. Hence, using $Z'(y)=A_*/(\wh
A(y) \wh W(y))$ the functional $\wh\calG$ from \eqref{eq:calG.A} is
transformed into $\calG$ via $\wh\calG(\alpha, U)= A_*\,\wh
\calG(\alpha/A_*, v)$.  The result of Proposition \ref{pr:Layer} now
follows from Proposition \ref{pr:MMM} via
\begin{align*}
\wh G(\alpha,u_0,u_1)&=\min\set{\wh\calG(\alpha,U)}{ U(0)=u_0,\ U(1)=u_1}\\
&=\min\set{ A_* \calG(\alpha/A_*,v)}{v(0)=U(0)/\wh W(0), \
  v(1)=U(1)/\wh W(1)}. 
\end{align*}
Thus, the asserted formula is established. 
\end{proof}

\subsection{Derivation of the potential $N$}
\label{ap:wh-N}

We consider the functional 
\[
\calN(v,\zeta)=\int_0^1 \Big(\frac{v'(z)^2}{2v(z)} - \frac 12
\zeta'(z)^2 v(z)\Big) \dd z
\]
for functions $v>0$. Hence, $\calN$ 
is convex in $v$ and concave in $\zeta$. 
We are interested in the inf-sup for given boundary values, namely
\[
N(\delta, v_0,v_1):= \inf\Bigset{ \sup\bigset{ \calN (v,\zeta)}{ 
   \zeta(1){-}\zeta(0)=\delta } }{ v(0)=v_0,\:v(1)=v_1,\: v>0 }. 
\]

The following result provides an explicit formula in terms of the dual
dissipation potential $\CCC^*$. It is based on the following surprising
relation between $\calN$ and $\calG$ from \eqref{eq:Galpha}:
\begin{equation}
  \label{eq:calM.calG}
\calM(\delta, v):= \max\bigset{ \calN (v,\zeta)}{ 
   \zeta(1){-}\zeta(0)=\delta }\ \overset{!!}{=}\ 
 \min\bigset{\calG(\alpha,v)-\alpha
   \delta}{\alpha \in \R}.
\end{equation}
The equality $\overset{!!}=$ can be checked by elementary
calculations, since in both cases we find 
\[
\calM(\delta,v)=\int_0^1 \frac{v'(z)^2}{2v(z)} \dd z -
\frac{\delta^2}2 \mafo{Harm}(v), \quad \text{where }\mafo{Harm}(v) =
\Big( \int_0^1 \frac{\rmd z}{v(z)}\Big)^{-1}.
\]
Using the strong link \eqref{eq:calM.calG} between $\calN$ and $\calG$ 
we show that $N$ can be calculated from $G$.

\begin{proposition}\label{pr:calN}
We have the relation
\[
N(\delta, v_0,v_1)  = 
 \sqrt{v_0v_1}\,\CCC^*(\log v_1{-}\log v_0) 
 - \sqrt{v_0v_1}\, \CCC^*(\delta).
\]
\end{proposition} 
\begin{proof}
  Using \eqref{eq:calM.calG} we want to show that $N$ is related
  to the Legendre transform
  $G^*(\delta,v_0,v_1):=\sup_{\alpha\in \R} \delta \alpha - G(\alpha,
  v_0,v_1)$ of $G$ from  \eqref{eq:Galpha}. 
For this we keep $\delta \in \R$ fixed. 
  
The functional $(\alpha,v)\mapsto \calG(\alpha,v) -\delta\alpha $ is
jointly convex, such that it can be minimized in any desired order of
$\alpha$ and $v$. Letting $V :=\set{v}{v>0,\;v(0)=v_0,\;v(1)=v_1}$ we
have 
\begin{align*}
  N(\delta,v_0,v_1)\ &=\  \inf_{v\in V} \calM(\delta,v)\ = \ 
   \inf_{v \in V}\Big( \inf_{\alpha\in \R} \calG(\alpha, v) {-} \delta \alpha\Big)   \\
  & =\  \inf_{\alpha\in \R}\Big( \inf_{v \in V} \calG(\alpha, v) {-}
  \delta \alpha\Big) = \inf_{\alpha\in \R}  G(\delta, v_0,v_1) -\delta\alpha \ = \
  - G^*(\delta, u_0,u_1) .
\end{align*}
Thus, evaluating $G^*$ with $G$ from \eqref{eq:A.MMM} explicitly gives
the desired result. 
\end{proof}
\medskip

\paragraph*{Acknowledgments.} M.L.\ was partially supported by the
Einstein Stiftung Berlin via the ECMath/\textsc{Matheon} project SE2.
A.M.\ was partially supported by 
DFG via project C5 within CRC 1114 (Scaling cascades in complex
systems) and by the ERC AdG.\,267802 \emph{AnaMultiScale}. 
M.R.\ was was partially supported by 
DFG via project C8 within CRC 1114 (Scaling cascades in complex
systems).

\footnotesize

\newcommand{\etalchar}[1]{$^{#1}$}
\def\cprime{$'$} \def\ocirc#1{\ifmmode\setbox0=\hbox{$#1$}\dimen0=\ht0
  \advance\dimen0 by1pt\rlap{\hbox to\wd0{\hss\raise\dimen0
  \hbox{\hskip.2em$\scriptscriptstyle\circ$}\hss}}#1\else {\accent"17 #1}\fi}

\end{document}